\documentclass{amsart}

\usepackage{mathtools}
\mathtoolsset{showonlyrefs=true}
\usepackage{amssymb}
\usepackage{bm}
\usepackage{graphicx}
\usepackage{subcaption}

\def\defeq{\mathrel{\mathop:}=} %:=

\newcommand{\Bm}{\mathrm{Bm}}
\newcommand{\Ca}{\mathrm{Ca}}
\newcommand{\relmiddle}[1]{\mathrel{}\middle#1\mathrel{}}

\begin{document}
    \title[Numerical scheme for Hele-Shaw problems by the MFS]{A simple numerical method for Hele-Shaw type problems by the method of fundamental solutions}
    \author[K.~Sakakibara]{Koya Sakakibara}
    \address[K.~Sakakibara]{Department of Applied Mathematics, Faculty of Science, Okayama University of Science, 1-1 Ridaicho, Kita-ku, Okayama-shi, Okayama 700-0005, Japan;
    RIKEN iTHEMS, 2-1 Hirosawa, Wako-shi, Saitama 351-0198, Japan}
    \email{ksakaki@xmath.ous.ac.jp}
    \author[Y.~Shimoji]{Yusauku Shimoji}
    \address[Y.~Shimoji]{Graduate School of Science and Technology, Meiji University, 1-1-1 Higashi-Mita, Tama-ku, Kawasaki-shi, Kanagawa 214-8571, Japan}
    \email{shimoji@meiji.ac.j}
    \author[S.~Yazaki]{Shigetoshi Yazaki}
    \address[S.~Yazaki]{Department of Mathematics, School of Science and Technology, Meiji University, 1-1-1 Higashi-Mita, Tama-ku, Kawasaki-shi, Kanagawa 214-8571, Japan}
    \email{syazaki@meiji.ac.jp}
    \begin{abstract}
        Hele-Shaw flows with time-dependent gaps create fingering patterns, and magnetic fluids in Hele-Shaw cells create intriguing patterns.
        We propose a simple numerical method for Hele-Shaw type problems by the method of fundamental solutions.
        The method of fundamental solutions is one of the mesh-free numerical solvers for potential problems, which provides a highly accurate approximate solution despite its simplicity.
        Moreover, combining with the asymptotic uniform distribution method, the numerical method satisfies the volume-preserving property.
        We use Amano's method to arrange the singular points in the method of fundamental solutions.
        We show several numerical results to exemplify the effectiveness of our numerical scheme.
    \end{abstract}
    \maketitle

    \section{Introduction}
    
    The Hele-Shaw problem describes the motion of a viscous fluid in a quasi-two-dimensional space, originating in a short paper by Henry Selby Hele-Shaw \cite{hele-shaw}, and was presented as a kind of experimental model for describing stream lines by using viscous fluids. 
    Let $\Omega(t)\subset\mathbb{R}^2$ be a bounded region occupied by a fluid, and $\Gamma(t)=\partial\Omega(t)$ be its boundary.
    Then, the Hele-Shaw problem is described as follows:
    \begin{align}
        \label{eq:HS-classic}
        \begin{dcases*}
            \triangle p(\cdot,t)=0&in $\Omega(t)$, $t\in(0,T)$,\\
            p(\cdot,t)=\sigma\kappa(\cdot,t)&on $\Gamma(t)$, $t\in(0,T)$,\\
            V(\cdot,t)=-h^2\nabla p(\cdot,t)\cdot\bm{N}(\cdot,t)&on $\Gamma(t)$, $t\in(0,T)$,
        \end{dcases*}
    \end{align}
    where $p(\cdot,t)$ is the pressure in $\Omega(t)$, $\sigma$ is the surface tension coefficient, $\kappa(\cdot,t)$ denotes the curvature, $V(\cdot,t)$ is the normal velocity, $\bm{N}(\cdot,t)$ is the unit outward normal vector of $\Gamma(t)$, and $h$ denotes the gap of the Hele-Shaw cell.
    The Laplace equation for pressure $p$, the first equation in \eqref{eq:HS-classic}, comes from the Darcy's law
    \begin{align}
        \bm{u}=-h^2\nabla p,\label{eq:Darcy}
    \end{align}
    and the incompressibility condition
    \begin{align}
        \nabla\cdot\bm{u}=0,\label{eq:incompressibility}
    \end{align}
    where $\bm{u}$ denotes the velocity field.
    The Dirichlet boundary condition, the second equation in \eqref{eq:HS-classic}, expresses the surface tension.
    For details of derivation, see, for example, the books \cite{lamb,gustafsson}.

    \begin{figure}[t]
        \centering
        \includegraphics[scale = 0.25]{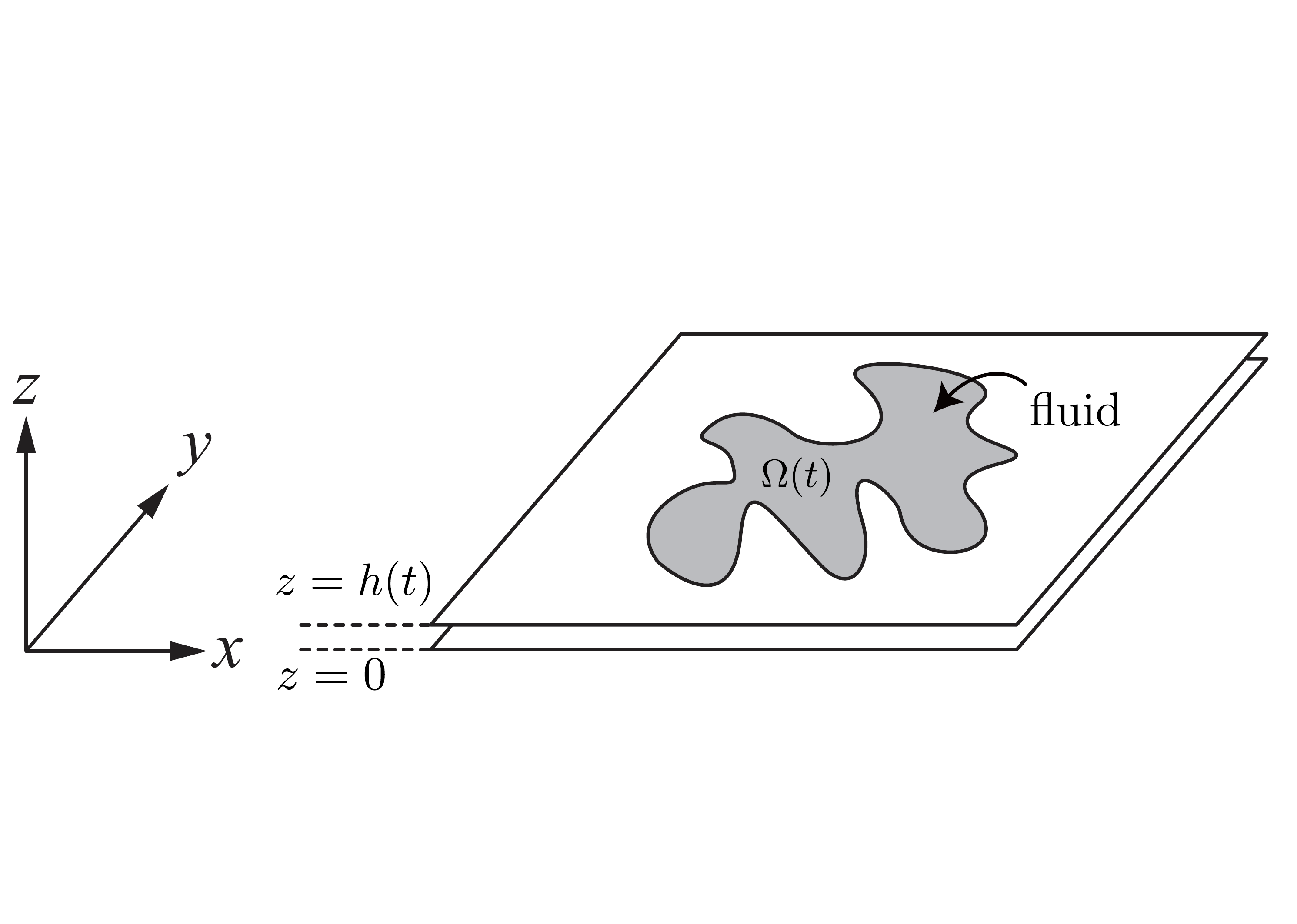}
        \caption{\label{Hele-Shaw_cell}Hele-Shaw cell}
    \end{figure}
    
    There are many variants of the Hele-Shaw problem, and in particular, they have been actively used for the experimental and mathematical study of fingering phenomena in various scientific fields \cite{saffman1958,howison,hou1994,shelley,tanveer2000,dockery}.
    Among them, Shelley et al. formulated the Hele-Shaw problem with a time-dependent gap in \cite{shelley} (see also Fig.\ref{Hele-Shaw_cell}):
    \begin{align}
        \label{eq:HS-TDG}
        \begin{dcases*}
            \triangle p(\cdot,t)=\frac{\dot{h}(t)}{h(t)^3}&in $\Omega(t)$, $t\in(0,T)$,\\
            p(\cdot,t)=\sigma\kappa(\cdot,t)&on $\Gamma(t)$, $t\in(0,T)$,\\
            V(\cdot,t)=-h(t)^2\nabla p(\cdot,t)\cdot\bm{N}(\cdot,t)&on $\Gamma(t)$, $t\in(0,T)$,
        \end{dcases*}
    \end{align}
    where $h(t)$ denotes the gap of the Hele-Shaw cell depending on time $t$.
    In the current situation, where the gap is time-dependent, the Darcy's law and the incompressibility condition are respectively given by
    \begin{align}
        \bm{u}=-h(t)^2\nabla p,
        \quad
        \nabla\cdot\bm{u}=-\frac{\dot{h}(t)}{h(t)}.
    \end{align}
    Therefore, the equation to be satisfied by the pressure is now the Poisson equation, not the Laplace equation.
    However, since a particular solution of the Poisson equation can be obtained concretely as
    \begin{align}
        p^*(\bm{x},t)
        =
        \frac{\dot{h}(t)}{4h(t)^3}|\bm{x}|^2,
    \end{align}
    by considering $p-p^*$ as $p$ again, \eqref{eq:HS-TDG} eventually leads to the following system:
    \begin{align}
        \label{eq:HS-TDG_dimensionless}
        \begin{dcases*}
            \triangle p(\cdot,t)=0&in $\Omega(t)$, $t\in(0,T)$\\
            p(\bm{x},t)=\sigma\kappa(\cdot,t)-\frac{\dot{h}(t)}{4h(t)^3}|\bm{x}|^2&for $\bm{x}\in\Gamma(t)$, $t\in(0,T)$,\\
            V(\bm{x},t)=-h(t)^2\left(\nabla p(\bm{x},t)+\frac{\dot{h}(t)}{2h(t)^3}\bm{x}\right)\cdot\bm{N}(\bm{x},t)&for $\bm{x}\in\Gamma(t)$, $t\in(0,T)$.
        \end{dcases*}
    \end{align}
    
    As a more advanced topic, fingering phenomena in magnetic fluids have also been investigated \cite{otto,rosensweig1983labyrinthine,elias1997macro,tatulchenkov}.
    The following is a summary of the mathematical model proposed in \cite{tatulchenkov}.
    Darcy's law and the incompressibility condition are respectively given by
    \begin{align}
        -\nabla p-\frac{12\eta}{h(t)^2}\bm{u}+\frac{2M_c}{h(t)}\nabla\varphi_m=\bm{0},
        \quad
        \nabla\cdot\bm{u}=-\frac{\dot{h}(t)}{h(t)},
    \end{align}
    where $\varphi_m=\varphi_m(\bm{x},t)$ is the magnetostatic field potential given by
    \begin{align}
        \varphi_m(\bm{x},t)=-M_c\int_{\Omega(t)}\left(\frac{1}{|\bm{x}-\bm{p}|}-\frac{1}{\sqrt{|\bm{x}-\bm{p}|^2+h(t)^2}}\right)\,\mathrm{d}\bm{p},
    \end{align}
    in which $M_c$ denotes the strength of the magnetization.
    In other words, we obtain the following system:
    \begin{align}
        \begin{dcases*}
            \triangle p(\cdot,t)=\frac{12\eta\dot{h}(t)}{h(t)^3}+\frac{2M_c}{h(t)}\triangle\varphi_m(\cdot,t)&in $\Omega(t)$, $t\in(0,T)$,\\
            p(\cdot,t)=\sigma\kappa(\cdot,t)&on $\Gamma(t)$, $t\in(0,T)$,\\
            V(\cdot,t)=-\frac{h(t)^2}{12\eta}\nabla\left(p-\frac{2M_c}{h}\varphi_m\right)(\cdot,t)\cdot\bm{N}(\cdot,t)&on $\Gamma(t)$, $t\in(0,T)$.
        \end{dcases*}
    \end{align}
    By constructing a particular solution of the Poisson equation in the same way as we derived \eqref{eq:HS-TDG_dimensionless}, and by performing the appropriate non-dimensionalization, we finally obtain the following model:
    \begin{align}
        \label{eq:HS_magnetic}
        \left\{
            \begin{array}{l}
                \triangle p(\cdot,t)=0\hspace{160pt} \text{in}\ \Omega(t),\ t\in(0,T),\\
                \displaystyle p(\cdot,t)=\kappa(\cdot,t)\left(\frac{R_0}{h_0}\right)^{1/3}\frac{\pi^{2/3}\Bm}{h_*(t)}\varphi(\cdot,t)-\left(\frac{R_0}{h_0}\right)^2\frac{\pi\Ca|\bm{x}|^2}{4h_*(t)^2}\\[2ex]
                \hspace{210pt} \text{on}\ \Gamma(t),\ t\in(0,T),\\[1ex]
                \displaystyle V(\bm{x},t)=-\frac{1}{\pi\Ca}\left(\frac{h_0}{R_0}\right)^2h_*(t)^2\nabla p(\bm{x},t)\cdot\bm{N}(\bm{x},t)-\frac{\dot{h}_*(t)}{2h_*(t)}\bm{x}\cdot\bm{N}(\bm{x},t)\\[2ex]
                \hspace{210pt} \text{for}\ \bm{x}\in\Gamma(t),\ t\in(0,T),
            \end{array}
        \right.
    \end{align}
    where $R_0$ is the radius of the initial drop, $h_0$ is the initial gap, $h_*(t)=\exp(t)$, $\varphi(\bm{x},t)=\varphi_m(\bm{x},t)/M_c$, $\Bm$ is the magnetic Bond number, and $\Ca$ is the capillary number.
    
    The magnetic fluid is one of the smart fluids.
    Because of the property which the smart fluids can be handled by controlling the electric and magnetic field, in other words, without touching directly, they are widely used in various fields, such as electric devices, mechanical engineering, and medical applications.
    In the 1960s, magnetic fluids were developed as part of space engineering development.
    R. E. Rosensweig is known as a pioneer of the investigation into magnetic fluids.
    His works (for instance, \cite{cowley1967,rosensweig1987}) are cited in many papers, even though they were contributed about half of a century ago.
    Rosensweig advocated the governing equation of magnetic fluid (see \cite{rosensweig1987} for more details):
    \begin{equation}
        \label{eq:governing magnetic fluid}
        \dfrac{\partial \bm{v}}{\partial t} + \left(\bm{v} \cdot \nabla \right) \bm{v} = -\nabla p^* + \mu_0 M \nabla H + \eta \triangle \bm{v} + \rho \bm{g}.
    \end{equation}
    This paper uses the magnetic fluid Hele-Shaw model \eqref{eq:HS_magnetic}, which can be derived from his equation \eqref{eq:governing magnetic fluid}.
    
    Magnetic fluid creates intriguing patterns. 
    For instance, many spikes are created in three-dimensional space by magnetic fluids in the magnetic field.
    The phenomena, called the spike phenomena, are well known as one of the intriguing phenomena caused by the collaboration of magnetic fluids and the magnetic field.
    On the other hand, the fluids create fingering and the so-called labyrinth patterns in a Hele-Shaw cell.
    Labyrinth patterns are mainly observed only in quasi-two-dimensional space.
    These two-dimensional and three-dimensional phenomena certainly have some common features; however, the difference of patterns shows differences between two-dimensional and three-dimensional problems.
    Therefore, proposing a simple and effective numerical method for magnetic fluid Hele-Shaw problems will significantly impact investigations into magnetic fluids.
    
    In order to solve these problems, various numerical methods have been proposed so far.
    The most successful method is based on the boundary integral method, originating in the celebrated paper \cite{hou1994} by T.Y.~Hou, J.S.~Lowengrub, and M.J.~Shelley.
    Because of the surface tension, explicit time discretization requires tiny time increments, and implicit discretization is difficult.
    By reformulating the problem differently, they overcame these difficulties and successfully computed the time evolution using a relatively large time increment.
    
    The original Hele-Shaw problem is known to have the following geometric variational structure.
    \begin{itemize}
        \item Curve-shortening property: the length $|\Gamma|(t)$ of $\Gamma(t)$ decays monotonically in time;
        \item Area-preserving property: the area $|\Omega|(t)$ of $\Omega(t)$ is constant in time;
        \item Barycenter-fixed property: the barycenter of $\Omega(t)$ does not move along time evolution.
    \end{itemize}
    In both the Hele-Shaw problem with time-dependent gap \eqref{eq:HS-TDG_dimensionless} and the Hele-Shaw problem for magnetic fluids \eqref{eq:HS_magnetic}, the volume-preserving property holds; that is, $|\Omega|(t)\cdot h(t)$ is constant in time.
    Therefore, it is natural to expect that these properties also hold in some sense in numerical methods.
    A method that preserves the mathematical structure of the target equation in a discrete sense is known as a structure-preserving numerical method, or geometric integrator, and has been the subject of many studies \cite{hairer}.
    
    Recently, the first and third authors of this paper proposed a structure-preserving numerical scheme for the Hele-Shaw problem in \cite{sakakibara2019}.
    The essential idea is to combine the method of fundamental solutions (MFS for short) with the uniform distribution method (UDM for short).
    The MFS is a numerical method for the Laplace equation that does not require any numerical integration and has the remarkable property that its approximation error decays exponentially for the number of sample points \cite{katsurada1988,1990katsurada,barnett2008,sakakibara2017}.
    The UDM is vital for stable numerical computation of moving boundary problems because it keeps the distance between adjacent nodes uniform with time evolution \cite{mikula2004direct,mikula2006evolution,sevcovic2011evolution,sevcovic2013gradient}.
    Furthermore, one of the main features of the numerical scheme proposed in \cite{sakakibara2019} is that, although it is in the sense of a semi-discrete problem, the area-preserving property is strictly achieved, and the curve-shortening and barycenter-fixed properties hold in an asymptotic manner.
    In this sense, the paper \cite{sakakibara2019} is the first structure-preserving numerical method to the Hele-Shaw problem.
    Other structure-preserving numerical schemes for moving boundary problems include, for example, \cite{kemmochi2017,sakakibara2021}.
    
    The purpose of this paper is to construct a numerical scheme satisfying the volume-preserving property for the Hele-Shaw problem with a time-dependent gap \eqref{eq:HS-TDG_dimensionless} and the Hele-Shaw problem with magnetic fluid \eqref{eq:HS_magnetic} by extending the contents of the paper \cite{sakakibara2019}.
    In order to solve problem \eqref{eq:HS_magnetic}, we need to compute the magnetostatic potential $\varphi$, i.e., the integral over $\Omega(t)$.
    In this paper, we adopt the Monte Carlo method for numerical integration.
    In addition, we improve the method adopted in \cite{sakakibara2019} to arrange the singular points that determine the accuracy of the approximate solution by the MFS.
    
    This paper is organized as follows.
    In section \ref{ss_scheme}, we construct a numerical scheme for problems \eqref{eq:HS-TDG_dimensionless} and \eqref{eq:HS_magnetic}.
    In particular, we show that the scheme satisfies the volume-preserving property as a semi-discrete scheme.
    Section \ref{result} reports numerical results for these problems and shows that the proposed method is much simpler than the previous ones but gives similar numerical results.
    Finally, in section \ref{conclusion}, we conclude this paper by presenting future research topics.
    
    \section{\label{ss_scheme}Numerical scheme}
\begin{figure}[tb]
        \centering
        \includegraphics[width=.7\hsize]{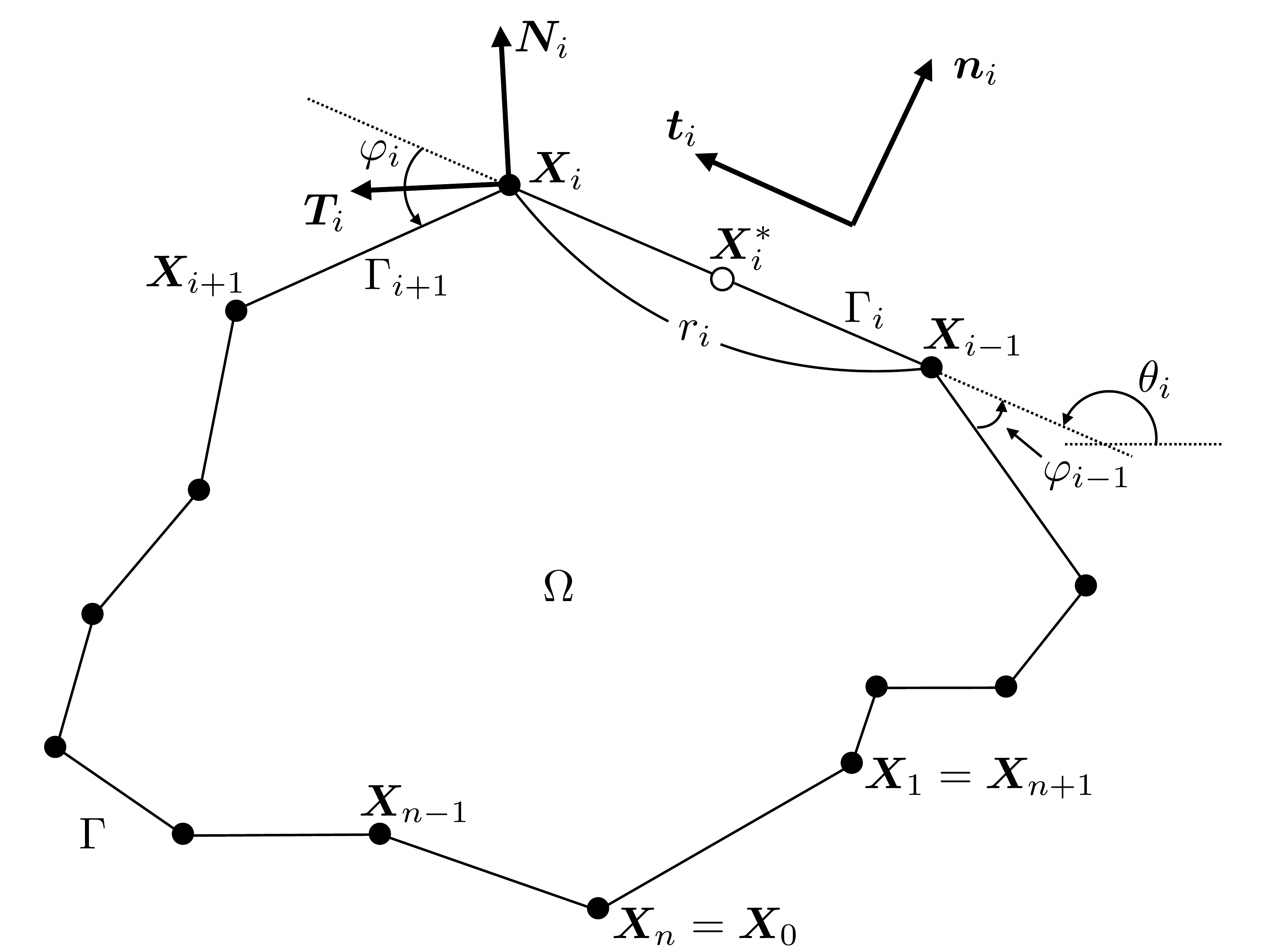}
        \caption{\label{poly_curve}Jordan polygonal curve}
    \end{figure}

    We approximate the boundary curve $\Gamma(t)$ by a polygonal curve; that is, let $\Gamma(t) = \bigcup_{i=1}^N \Gamma_i(t)$ be an $N$-sided polygonal Jordan curve, where $\Gamma_i(t)$ is the $i$-th edge of $\Gamma(t)$ defined by
    \begin{equation}
        \Gamma_i(t) = (\bm{X}_{i-1}(t), \bm{X}_i(t)) \defeq \{\lambda \bm{X}_{i-1}(t) + (1 - \lambda)\bm{X}_i(t) \mid \lambda \in (0, 1)\},
    \end{equation}
    where $\bm{X}_i(t)$ is the $i$-th vertex of $\Gamma(t)$ ($i = 1, 2, \ldots, N$.). These $N$ vertices are labeled in the anti-clockwise order throughout this paper, and we adopt the periodic notation $\bm{X}_0(t)=\bm{X}_N(t)$ and $\bm{X}_{N+1}(t)=\bm{X}_1(t)$ (see Fig.(\ref{poly_curve})).
    The motion of the polygonal curve is then described by the following evolution law, which is a system of ordinary differential equations:
    \begin{equation}
        \label{dis_evo_law}
        \dot{\bm{X}}_i(t) = V_i(t) \bm{N}_i(t) + W_i \bm{T}_i(t),
        \quad
        i=1,2,\ldots,N,\ t\in(0,T),
    \end{equation}
    where $\bm{N}_i$ and $\bm{T}_i$ are the unit outward normal and unit tangent vectors at $\bm{X}_i$, and $V_i$ and $W_i$ denote the normal and tangential velocities at $\bm{X}_i$, respectively. 
    
    We here summarize our numerical scheme.
    The details of each step are described in the subsequent subsections.
    
    For a given polygonal curve $\Gamma^{(m)}$ at the $m$-th time, we compute the
    right-hand side of the evolution law \eqref{dis_evo_law} as in the following three steps.
    \begin{description}
        \item[\textbf{Step 1:}]
        Define the unit tangent vector $\bm{T}_i^{(m)}$ at $\bm{X}_i^{(m)}$ as the unit vector pointing in the exact middle direction of the unit tangent vectors defined on the two adjacent edges.
        Then, by rotating it by $-\pi/2$, we define the unit outward normal vector $\bm{N}_i^{(m)}$ at $\bm{X}_i^{(m)}$.
        We also define the discrete curvature $\kappa_i^{(m)}$ on each edge as the first variation of the perimeter.
        \item[\textbf{Step 2:}]
        The MFS solves the Dirichlet problem, \eqref{eq:HS-TDG_dimensionless} or \eqref{eq:HS_magnetic}, in the polygonal region $\Omega^{(m)}$.
        For example, \eqref{eq:HS-TDG_dimensionless} can be formulated as in the following using the discrete curvatures $\{\kappa_i^{(m)}\}_{i=1}^N$ defined in Step 1:
        \begin{align}
            \begin{dcases*}
                \triangle p(\cdot,t^{(m)})=0&in $\Omega^{(m)}$,\\
                p(\bm{x},t^{(m)})=\sigma\kappa_i^{(m)}-\frac{\dot{h}(t^{(m)})}{4h(t^{(m)})^3}|\bm{x}|^2&for $\bm{x}\in\Gamma_i^{(m)}$, $i=1,2,\ldots,N$.
            \end{dcases*}
        \end{align}
        Using the approximate solution $P^{(m)}$ by the MFS and the unit outward normal vector $\bm{N}_i^{(m)}$ at $\bm{X}_i^{(m)}$, the normal velocity $V_i^{(m)}$ can be computed.
        The same procedure also works for problem \eqref{eq:HS_magnetic}.
        \item[\textbf{Step 3:}]
        Compute the tangential velocities $\{W_i^{(m)}\}_{i=1}^N$ using the uniform distribution method.
    \end{description}
    Since the right-hand side of equation \eqref{dis_evo_law} can be computed by the above steps, the evolution law \eqref{dis_evo_law} is computed by the fourth-order Runge--Kutta method to obtain the polygonal curve $\Gamma^{(m+1)}$ at the next time $t^{(m+1)}$.
    
    In the following, superscripts indicating time steps are omitted so long as there is no risk of confusion.
    
    \subsection{Step 1: Spatial discretization}
    
    \subsubsection{Definitions of $\bm{T}_i$ and $\bm{N}_i$}
    Let $r_i$ be the length of $\Gamma_i$ and let $\bm{t}_i$ and $\bm{n}_i$ be the unit tangent vector and unit outward normal vector defined naturally on $\Gamma_i$.
    
    Unit tangent vectors and unit outward normal vectors at vertices are not well-defined, and there are various possible ways of defining them.
    In this paper, we adopt the following definitions.
    Let $\varphi_i$ be a signed exterior angle at $\bm{X}_i$; that is, the angle between two adjacent edges $\Gamma_i$ and $\Gamma_{i+1}$, satisfying $\bm{t}_i\cdot\bm{t}_{i+1}=\cos\varphi_i$.
    Using this $\varphi_i$, we define the unit tangent vector $\bm{T}_i$ as follows:
    \begin{align}
        \bm{T}_i
        =
        \frac{\bm{t}_i+\bm{t}_{i+1}}{|\bm{t}_i+\bm{t}_{i+1}|}
        =
        \frac{\bm{t}_i+\bm{t}_{i+1}}{2\mathsf{cos}_i},
        \quad
        i=1,2,\ldots,N,
        \label{eq:def_T}
    \end{align}
    where $\mathsf{cos}_i\coloneqq\cos(\varphi_i/2)$.
    The unit outward normal vector $\bm{N}_i$ is then defined by $\bm{N}_i=-\bm{T}_i^\bot$.
    
    \subsubsection{Definition of discrete curvature}
    The curvature is usually defined as satisfying the Frenet formula.
    However, since each side of a polygonal curve is a line segment, the usual definition of curvature will lead to zero curvature almost everywhere, which is not suitable to study the motion of polygonal curves.
    On the other hand, the curvature can be interpreted as the first variation of the length functional.
    This paper applies this view to polygonal curves and derives an appropriate expression for discrete curvature.
    
    The perimeter $L$ of a polygonal curve $\Gamma$ is given by
    \begin{align}
        L=\sum_{i=1}^Nr_i.
    \end{align}
    According to the evolution law \eqref{dis_evo_law}, the time derivative of $L$ can be computed as follows:
    \begin{align}
        \dot{L}=2\sum_{i=1}^NV_i\mathsf{sin}_i,
        \label{eq:dot_L_pre}
    \end{align}
    where $\mathsf{sin}_i=\sin(\varphi_i/2)$.
    Let us assume that the following relation holds between the normal velocities $\{V_i\}_{i=1}^N$ at vertices and the normal velocities $\{v_i\}_{i=1}^N$ on edges:
    \begin{align}
        V_i=\frac{v_i+v_{i+1}}{2\mathsf{cos}_i},
        \quad
        i=1,2,\ldots,N,
        \label{eq:relation_V_v}
    \end{align}
    which is an analogy of \eqref{eq:def_T}.
    Then, the equation \eqref{eq:dot_L_pre} for the time derivative of the perimeter is
    \begin{align}
        \dot{L}=\sum_{i=1}^N\kappa_iv_ir_i,
        \quad
        \kappa_i=\frac{\mathsf{tan}_i+\mathsf{tan}_{i-1}}{r_i},
        \quad
        i=1,2,\ldots,N,
    \end{align}
    where $\mathsf{tan}_i=\mathsf{sin}_i/\mathsf{cos}_i$.
    We define $\kappa_i$ to be the discrete curvature of the $i$-th edge $\Gamma_i$.
    
    \subsection{Step 2: Computation of normal velocities by the MFS}
    
    Normal velocities are computed using the MFS, a mesh-free numerical solver for potential problems.
    For example, in the Hele-Shaw problem with a time-dependent gap \eqref{eq:HS-TDG_dimensionless}, the boundary is approximated by a polygonal curve, and the discrete curvature gives the curvature on edge, so the potential problem to be solved at time $t=t^{(m)}$ is
    \begin{align}
        \begin{dcases*}
            \triangle p(\cdot,t^{(m)})=0&in $\Omega^{(m)}$,\\
            p(\bm{x},t^{(m)})=\sigma\kappa_i^{(m)}-\frac{\dot{h}(t^{(m)})}{4h(t^{(m)})^3}|\bm{x}|^2&for $\bm{x}\in\Gamma_i^{(m)}$, $i=1,2,\ldots,N$.
        \end{dcases*}
    \end{align}
    Following the paper \cite{sakakibara2019}, the approximate solution $P$ based on the MFS is given by
    \begin{align}
        P(\bm{x})=Q_0+\sum_{j=1}^NQ_jE_j(\bm{x}),
        \quad
        E_j(\bm{x})=E(\bm{x}-\bm{y}_j)-E(\bm{x}-\bm{z}_j),
    \end{align}
    where $E(\bm{x})=(2\pi)^{-1}\log|\bm{x}|$ is the fundamental solution of the Laplace operator.
    The points $\{\bm{y}_j\}_{j=1}^N$ and $\{\bm{z}_j\}_{j=1}^N$ are singular points and dummy points, respectively, chosen ``suitably'' from outside the region $\Omega^{(m)}$.
    They play a role in controlling the accuracy of the MFS approximate solution.
    The specific arrangement of these points will be discussed in subsection \ref{subsec:Amano}.
    Since $P$ satisfies the Laplace equation exactly, the linear combination coefficients $\{Q_j\}_{j=0}^N$ are determined by an approximate treatment of the Dirichlet boundary conditions and the volume-preserving property.
    In this study, we adopt the midpoint $\bm{X}_i^*=(\bm{X}_i+\bm{X}_{i-1})/2$ of each edge as a collocation point and approximate the Dirichlet boundary condition using the collocation method.
    In other words, $\{Q_j\}_{j=0}^N$ are given as a solution to the following linear system
    \begin{align}
        P(\bm{X}_i^*)
        =
        \sigma\kappa_i^{(m)}-\frac{\dot{h}(t^{(m)})}{4h(t^{(m)})^3}|\bm{X}_i^*|^2,
        \quad
        i=1,2,\ldots,N,
    \end{align}
    with a constraint
    \begin{align}
        \sum_{j=1}^NQ_jH_j=0,
        \quad
        H_j=\sum_{i=1}^N\nabla E_j(\bm{X}_i^*)\cdot\bm{n}_ir_i,
        \quad
        j=1,2,\ldots,N.
        \label{eq:constraint}
    \end{align}
    
    By differentiating the approximate solution $P$ in the normal direction $\bm{n}_i$, the normal velocity $v_i$ on the $i$-th edge can be computed by
    \begin{align}
        v_i=-h(t^{(m)})^2\left(\nabla P(\bm{X}_i^*)+\frac{\dot{h}(t^{(m)})}{2h(t^{(m)})}\bm{X}_i^*\right)\cdot\bm{n}_i,
        \quad
        i=1,2,\ldots,N.
        \label{eq:normal_velocity_edge}
    \end{align}
    The great advantage of using the MFS is that the gradient $\nabla P$ can be computed analytically.
    In other words, there is no need to perform numerical differentiation.
    Finally, the normal velocities $\{V_i\}_{i=1}^N$ at vertices can be obtained by using the relation \eqref{eq:relation_V_v}.
    
    At the end of this subsection, we clarify what constraint \eqref{eq:constraint} means.
    To this end, we state a little more precisely the volume-preserving property.
    Since the volume $V(t)$ of the fluid region $\Omega(t)\times(0,h(t))$ is given by $V(t)=A(t)h(t)$, where $A(t)=|\Omega|(t)$, the volume-preserving property $\dot{V}=0$ can be rephrased as the following equation for the time evolution of the area:
    \begin{align}
        \dot{A}(t)=-\frac{\dot{h}(t)}{h(t)}A(t).
        \label{eq:volume-preserving_area}
    \end{align}
    According to Proposition 4 in the paper \cite{sakakibara2019}, the time evolution of the area under the evolution law \eqref{dis_evo_law} is given by
    \begin{align}
        \dot{A}=\sum_{i=1}^Nv_ir_i+\mathrm{err}_A,
        \quad
        \mathrm{err}_A=\sum_{i=1}^N\left(W_i\mathsf{sin}_i-\frac{v_{i+1}-v_i}{2}\right)\frac{r_{i+1}-r_i}{2}.
        \label{eq:error_term}
    \end{align}
    If the error term $\mathrm{err}_A$ is equal to $0$, then the normal velocity formula \eqref{eq:normal_velocity_edge} gives the time evolution of the area as
    \begin{align}
        \dot{A}(t)
        =
        -h(t)^2\sum_{j=1}^NQ_jH_j-\frac{\dot{h}(t)}{h(t)}A(t).
    \end{align}
    Therefore, because of \eqref{eq:volume-preserving_area}, we obtain constraint \eqref{eq:constraint}.
    
    For the above argument to be consistent, the error term $\mathrm{err}_A$ must be equal to $0$.
    For this purpose, we adopt the uniform distribution method to compute the tangential velocity.
    Details are given in subsection \ref{subsec:UDM}.
    
    The same consideration can be applied to the Hele-Shaw problem for magnetic fluids \eqref{eq:HS_magnetic} to obtain an approximate solution satisfying the volume-preserving property.
    
    \subsection{Step 3: Computation of tangential velocities by the UDM}
    \label{subsec:UDM}
    As mentioned in the previous subsection, the error term must be $0$ to achieve the volume-preserving property.
    Given the representation \eqref{eq:error_term} of the error term, one possibility is to apply the UDM.
    It is a method of moving vertices in the tangential direction so that the length of each edge of the polygon is equal, i.e., $r_i=L/N$ ($i=1,2,\ldots,N$), and is known to provide stable numerical computations.
    In this study, we adopt the asymptotic UDM, which achieves uniform distribution in an asymptotic sense.
    This method is robust to the effects of numerical errors such as rounding errors, and the tangential velocities are given as linear expressions of normal velocities.
    The asymptotic UDM requires that $r_i(t)\rightarrow L(t)/N$ as $t\to T_{\max}$, not $r_i(t)=L(t)/N$, holds, where $T_{\max}$ represents a final computation time.
    
    We assume that the following relations hold to derive an expression for the tangential velocity by the asymptotic UDM:
    \begin{equation}
        r_i - \dfrac{L}{N} = \eta_i e^{-\mu(t)},\quad i=1,2,\ldots,N,
        \label{eq:UDM_1}
    \end{equation}
    where $\{\eta_i\}_{i=1}^N$ are parameters chosen so that their mean is equal to $0$.
    The function $\mu$ controls the strength of the effect of the uniform arrangement, and a candidate satisfies $\mu(t)\rightarrow\infty$ as $t\to T_{\max}$.
    Differentiating both sides of equation \eqref{eq:UDM_1} with respect to $t$ yields
    \begin{equation}
        \dot{r}_i = \dfrac{\dot{L}}{N} + \left(\dfrac{L}{N} - r_i\right)\omega(t), \quad i=1,2,\ldots,N,
        \label{eq:dot_r}
    \end{equation}
    where $\omega(t)=\dot{\mu}(t)$.
    On the other hand, if we compute the time derivative of the length of the edge using the evolution law \eqref{dis_evo_law}, we obtain
    \begin{align}
        \dot{r}_i=V_i\mathsf{sin}_i+V_{i-1}\mathsf{sin}_{i-1}+W_i\mathsf{cos}_i-W_{i-1}\mathsf{cos}_{i-1},
        \quad
        i=1,2,\ldots,N.
        \label{eq:UDM_2}
    \end{align}
    Combining \eqref{eq:dot_r} and \eqref{eq:UDM_2}, we obtain recurrence relations describing the tangential velocity:
    \begin{equation*}
        W_i \mathsf{cos}_i - W_{i-1} \mathsf{cos}_{i-1} = \frac{\dot{L}}{N}+\left(\frac{L}{N}-r_i\right)\omega(t) - V_i \mathsf{sin}_i - V_{i-1} \mathsf{sin}_{i-1}
    \end{equation*}
    for $i = 2,3, \ldots, N$. 
    Imposing the zero-average condition $\sum_{i=1}^ N W_i = 0$, we obtain $N$ linearly independent equations.
    This system can be solved analytically, resulting in the following expression for the tangential velocity in the asymptotic UDM:
    \begin{align}
        W_i=\frac{\Psi_i+W_1\mathsf{cos}_1}{\mathsf{cos}_i},
        \quad
        i=2,3,\ldots,N,
        \quad
        W_1=-\frac{\sum_{i=2}^N\Psi_i/\mathsf{cos}_i}{\mathsf{cos}_1\sum_{j=1}^N\mathsf{cos}_j^{-1}},
        \label{eq:UDM}
    \end{align}
    where
    \begin{align}
        &\Psi_i=\sum_{l=2}^i\psi_l,
        \quad
        \psi_i=-V_i\mathsf{sin}_i-V_{i-1}\mathsf{sin}_{i-1}+\frac{\dot{L}}{N}+\left(\frac{L}{N}-r_i\right)\omega
    \end{align}
    for $i=2,3,\ldots,N$.
    Note that if the initial arrangement of vertices is uniform, then the asymptotic UDM theoretically keeps uniformity during time evolution.
    
    There is another way to eliminate $\mathrm{err}_A$.
    If $W_i$ is defined by
    \begin{align}
        W_i=\frac{v_{i+1}-v_i}{2\mathsf{sin}_i},
        \quad
        i=1,2,\ldots,N,
    \end{align}
    then the error term $\mathrm{err}_A$ also vanishes.
    However, in this case, the orientation of the edges of the polygonal curve does not change in time. 
    Namely, the time evolution of the polygon is restricted to the crystalline motion.
    This tangential velocity is useful for theoretical analysis, but it is not suitable for numerical computation of problems with large deformation of the domain, such as the Hele-Shaw problem.
    Therefore, in this study, we adopt the asymptotic UDM given by equation \eqref{eq:UDM}.
    
    \section{\label{result}Numerical result}
    In this section, we show some results of our numerical computations for \eqref{eq:HS-TDG_dimensionless} and \eqref{eq:HS_magnetic}.
    
    \subsection{Arrangement of points by the modified Amano's method}
    \label{subsec:Amano}
    In solving problems \eqref{eq:HS-TDG_dimensionless} or \eqref{eq:HS_magnetic} by the MFS, singular points and dummy points must be appropriately placed.
    In the paper \cite{sakakibara2019}, the singular points are placed at a fixed distance in the unit outward normal direction from the midpoints of the edges:
    \begin{align}
        \bm{y}_j=\bm{X}_j^*+d\bm{n}_j,
        \quad
        j=1,2,\ldots,N.
        \label{eq:singular_old}
    \end{align}
    However, there is no mathematical background to this point placement, and there may be better ways to arrange them.
    
\begin{figure}[tb]
        \centering
        \includegraphics[width=.7\hsize]{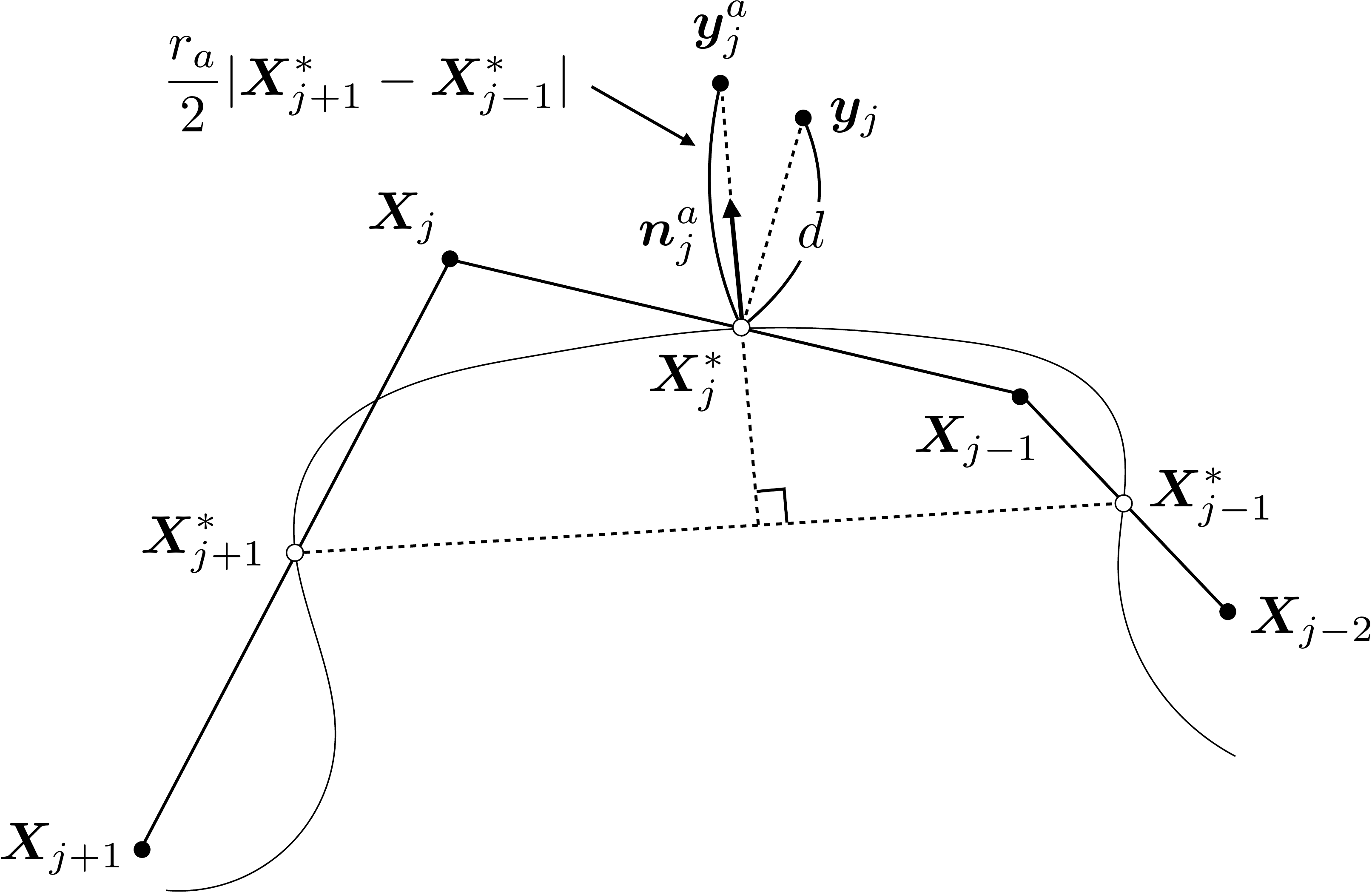}
        \caption{Arrangement of the singular point.
        $\bm{y}_j$ corresponds to the original arrangement \eqref{eq:singular_old} adopted in \cite{sakakibara2019}, while $\bm{y}_j^a$ does to the modified Amano's arrangement \eqref{eq:singular_Amano}.}
        \label{fig:singular}
    \end{figure}

    In this study, we employ a modified Amano's method proposed in numerical computation of conformal mappings (see for instance \cite{amano2012}).
    Namely, the singular points are placed according to the following relations:
    \begin{align}
        \bm{y}_j^a=\bm{X}_j^*+\frac{r_a}{2}|\bm{X}_{j+1}^*-\bm{X}_{j-1}^*|\bm{n}_j^a,
        \label{eq:singular_Amano}
    \end{align}
    where $\bm{n}_j^{a}$ is a modified unit outward normal vector determined by the positions of the three points $\bm{X}_{j-1}^*$, $\bm{X}_j^*$, and $\bm{X}_{j+1}^*$ and is defined through
    \begin{align}
        \bm{n}_j^a
        =
        -\frac{(\bm{X}_{j+1}^*-\bm{X}_{j-1}^*)^\bot}{|\bm{X}_{j+1}^*-\bm{X}_{j-1}^*|}.
    \end{align}
    The geometric picture of the arrangement of singular points is shown in Fig. \ref{fig:singular}.
    
    Amano's method can be though of as corresponding to a discrete version of point arrangement using a conformal mapping, which is used in the theoretical analysis of the MFS.
    For example, in Fig. \ref{fig:singular}, the thin line expresses the original smooth curve.
    In this case, $\bm{n}_j^a$ is considered an approximation of the unit outward normal vector at $\bm{X}_i^*$.
    This perspective is correct, and it can be mathematically proved that Amano's method approximates the method using conformal mappings.
    For more details, refer to the paper \cite{sakakibara2020}.
    Hence, we use the modified Amano's arrangement \eqref{eq:singular_Amano} to place the singular points in this paper.
    
    While the singular points have a significant impact on the accuracy of the approximate solution, the dummy points are only introduced for the invariance of the approximate solution and have little effect on the accuracy of the approximate solution.
    Therefore, in this paper, we adopt the following simple arrangement:
    \begin{align}
        \bm{z}_j
        =
        1000(\cos\alpha_i,\sin\alpha_i)^\top,
        \quad
        \alpha_i=\frac{2\pi i}{N},
        \quad
        i=1,2,\ldots,N.
    \end{align}
    
    \subsection{Hele-Shaw problem with a time-dependent gap}
    
    We show the results of numerical computations for \eqref{eq:HS-TDG_dimensionless}.
    The initial curve $\Gamma(0) : [0, 1] \ni u \mapsto (x_1(u), x_2(u))^\top \in \mathbb{R}^2$ is given by
    \begin{align*}
        &x_1(u) = r\cos(2\pi u), \quad x_2 = r\sin(2\pi u),\\
        &r = R_0 + 0.02(\cos(6\pi u) + \sin(14\pi u) + \cos(30\pi u) + \sin(50\pi u)).
    \end{align*}
    We set the parameters as in Table~\ref{tab_TDG}.
    \begin{table}[tb]
        % \centering
        \caption{\label{tab_TDG}parameters for the Hele-Shaw problem with a time-dependent gap \eqref{eq:HS-TDG_dimensionless}}
        \begin{tabular}{l|l}
        \hline \hline
        $N = 300$ & the number of vertices\\
        $\sigma = 2.0 \times 10^{-4}$ & the surface tension coefficient\\
        $\Delta t = 1 / (10N^2)$ & the time-mesh size\\
        $r_a = 1$ & the parameter in the Amano's method \eqref{eq:singular_Amano}\\
        $R_0 = 1$ & the initial reference radius\\
        \hline
        \end{tabular}
    \end{table}
    
\begin{figure}[tb]
        \begin{tabular}{cccc}
            \begin{minipage}{0.25\linewidth}
                \centering
                \includegraphics[scale = 0.32]{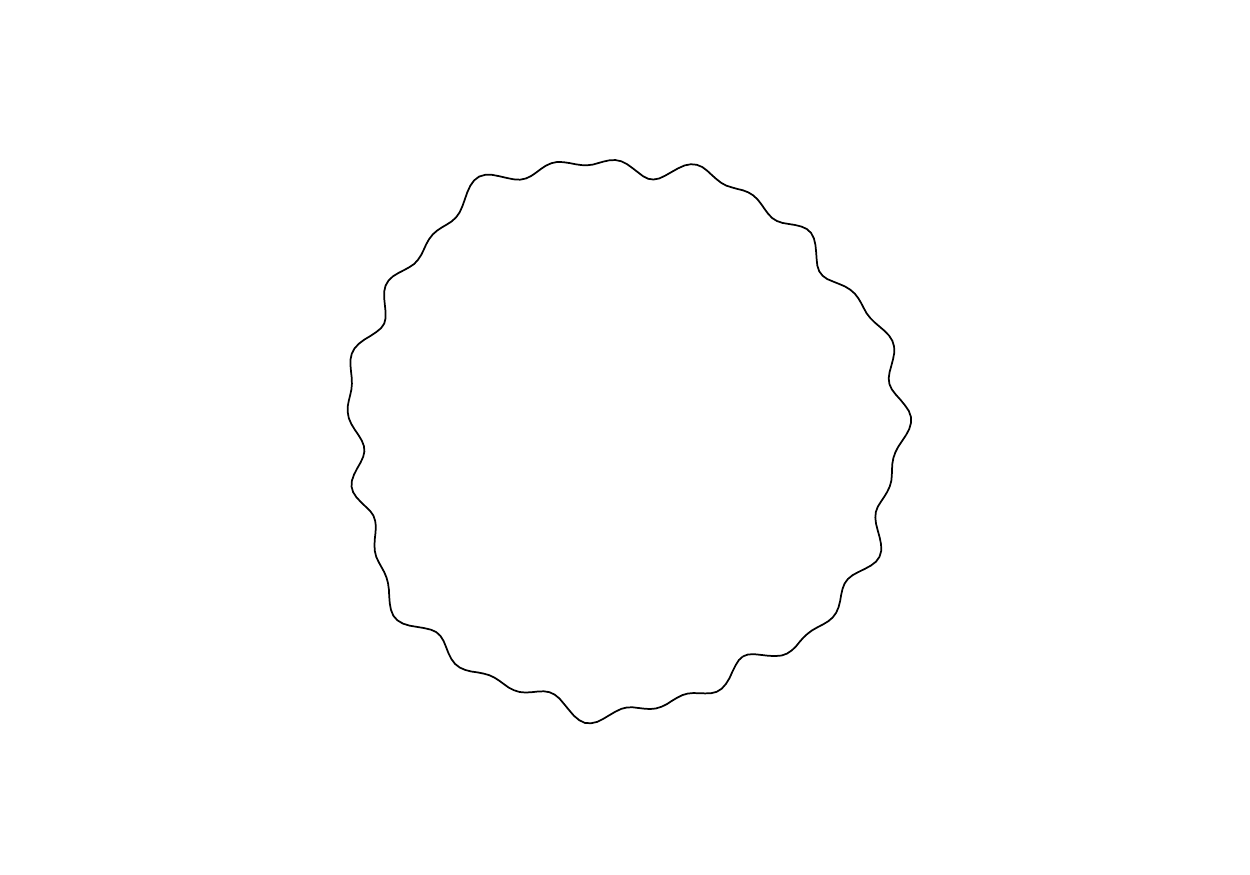}
                \subcaption{$t=0.0$}
            \end{minipage}
            \begin{minipage}{0.25\linewidth}
                \centering
                \includegraphics[scale = 0.32]{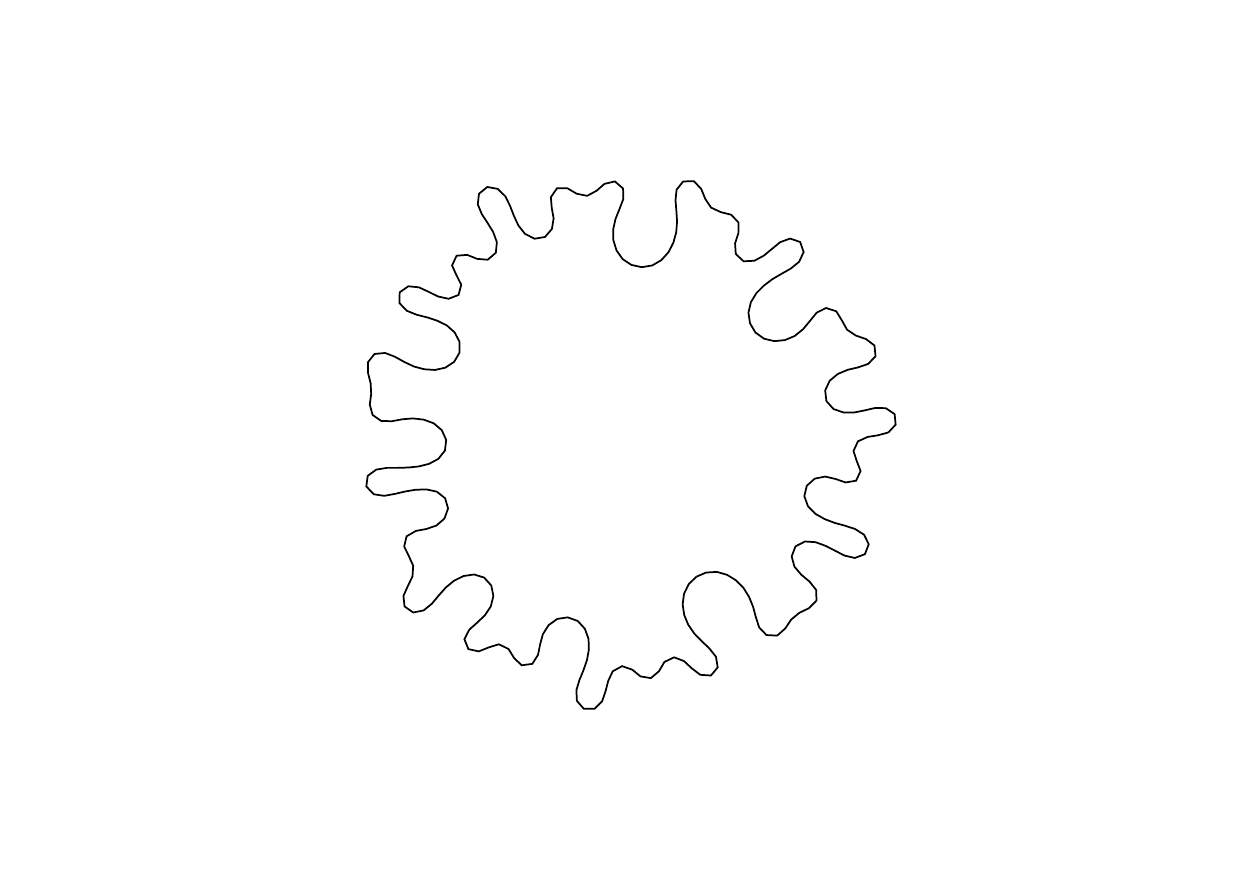}
                \subcaption{$t=0.36$}
            \end{minipage}
            \begin{minipage}{0.25\linewidth}
                \centering
                \includegraphics[scale = 0.32]{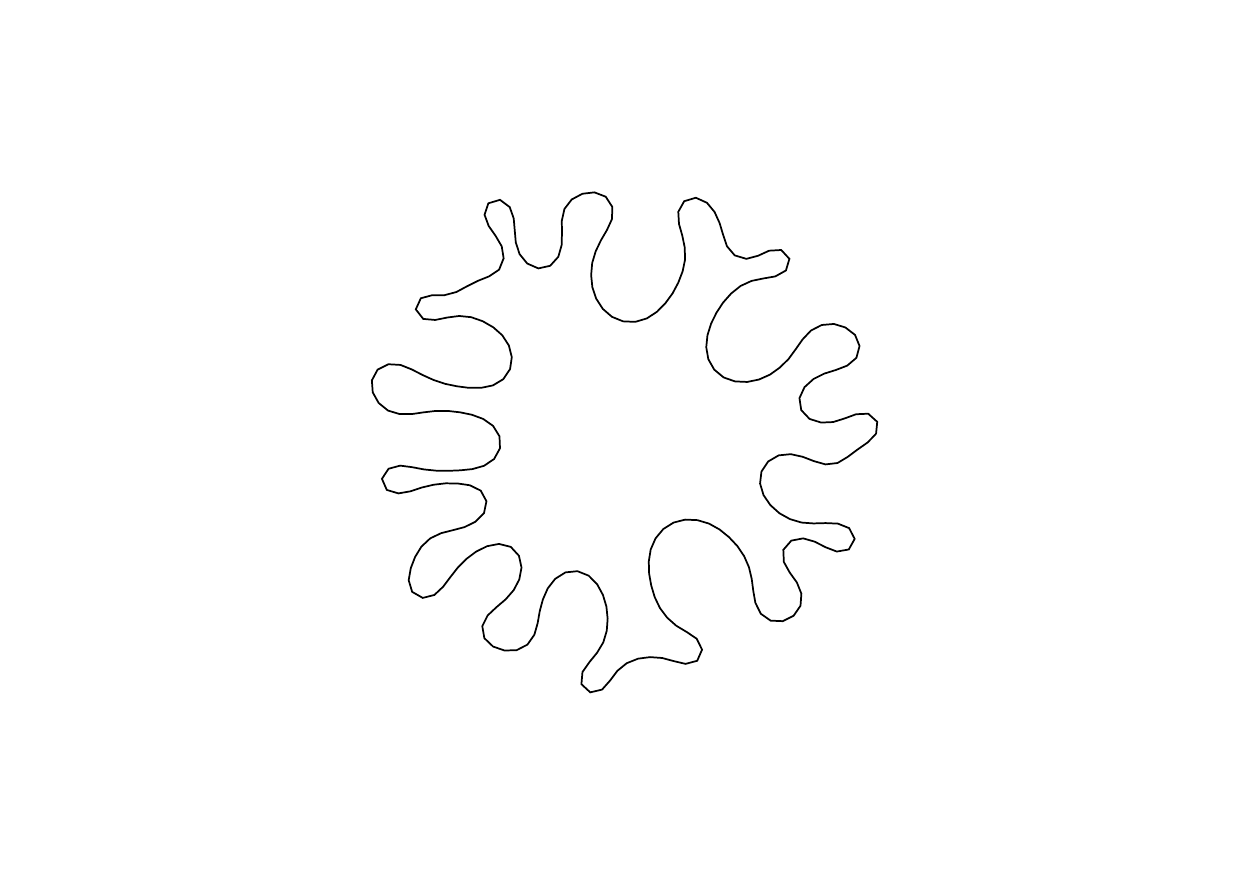}
                \subcaption{$t=0.73$}
            \end{minipage}
            \begin{minipage}{0.25\linewidth}
                \centering
                \includegraphics[scale = 0.32]{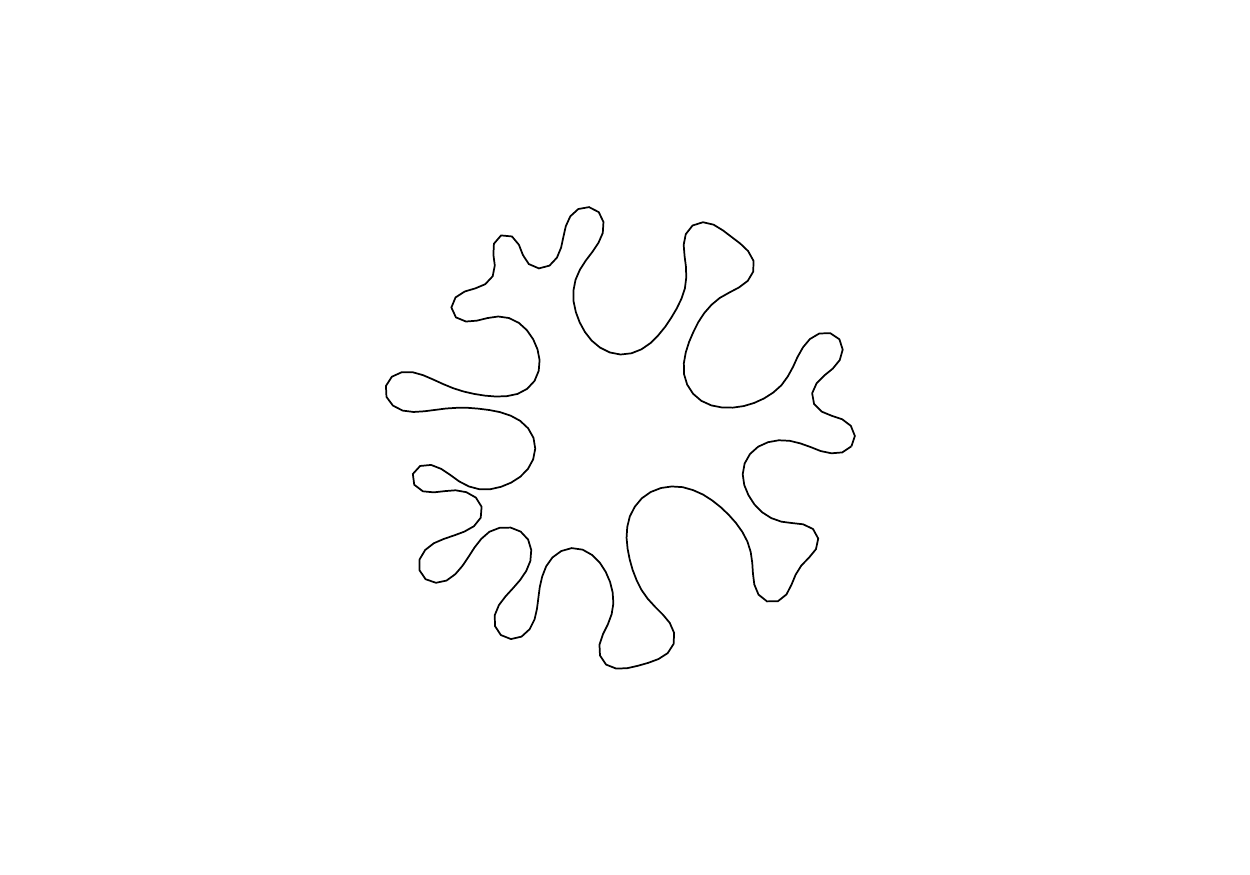}
                \subcaption{$t=1.10$}
            \end{minipage}
            \\
            \begin{minipage}{0.25\linewidth}
                \centering
                \includegraphics[scale = 0.32]{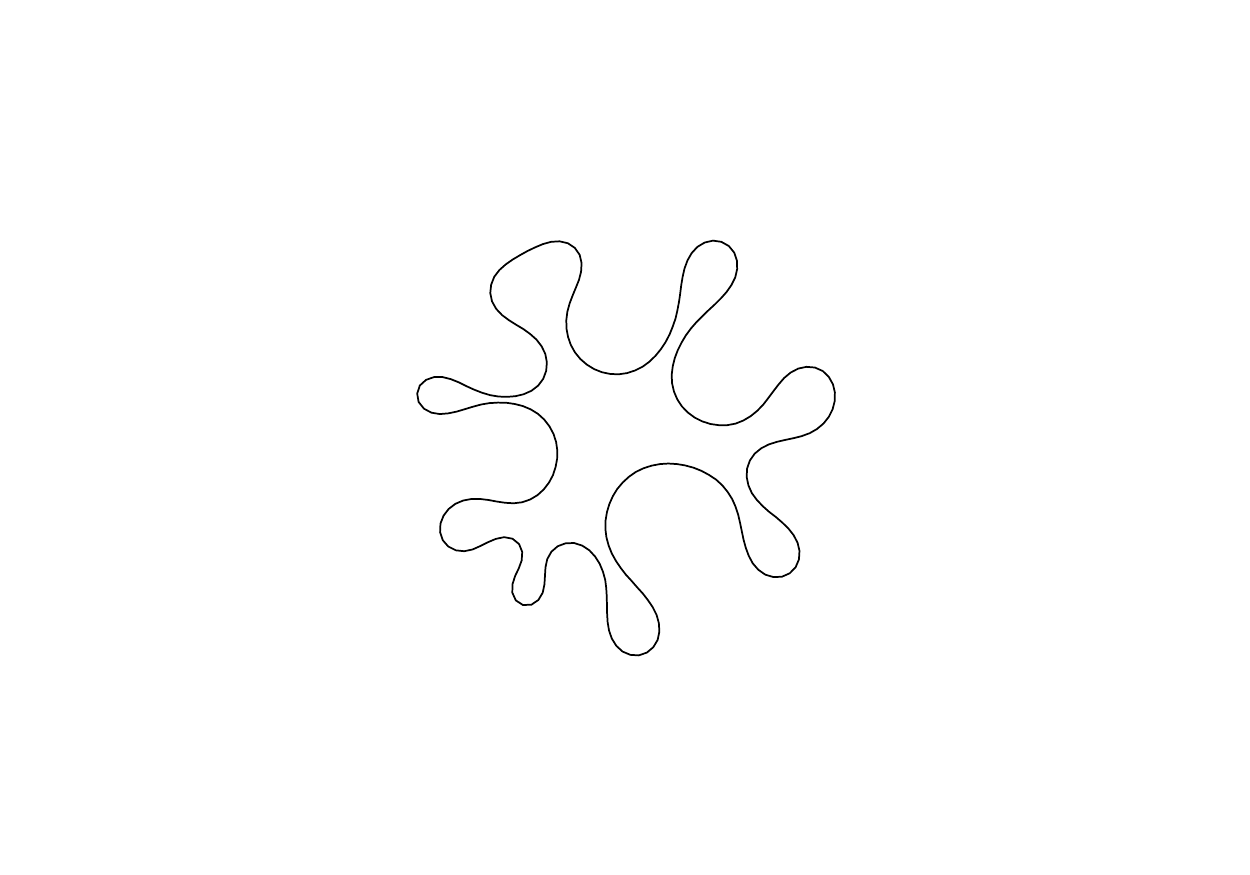}
                \subcaption{$t=1.47$}
            \end{minipage}
            \begin{minipage}{0.25\linewidth}
                \centering
                \includegraphics[scale = 0.32]{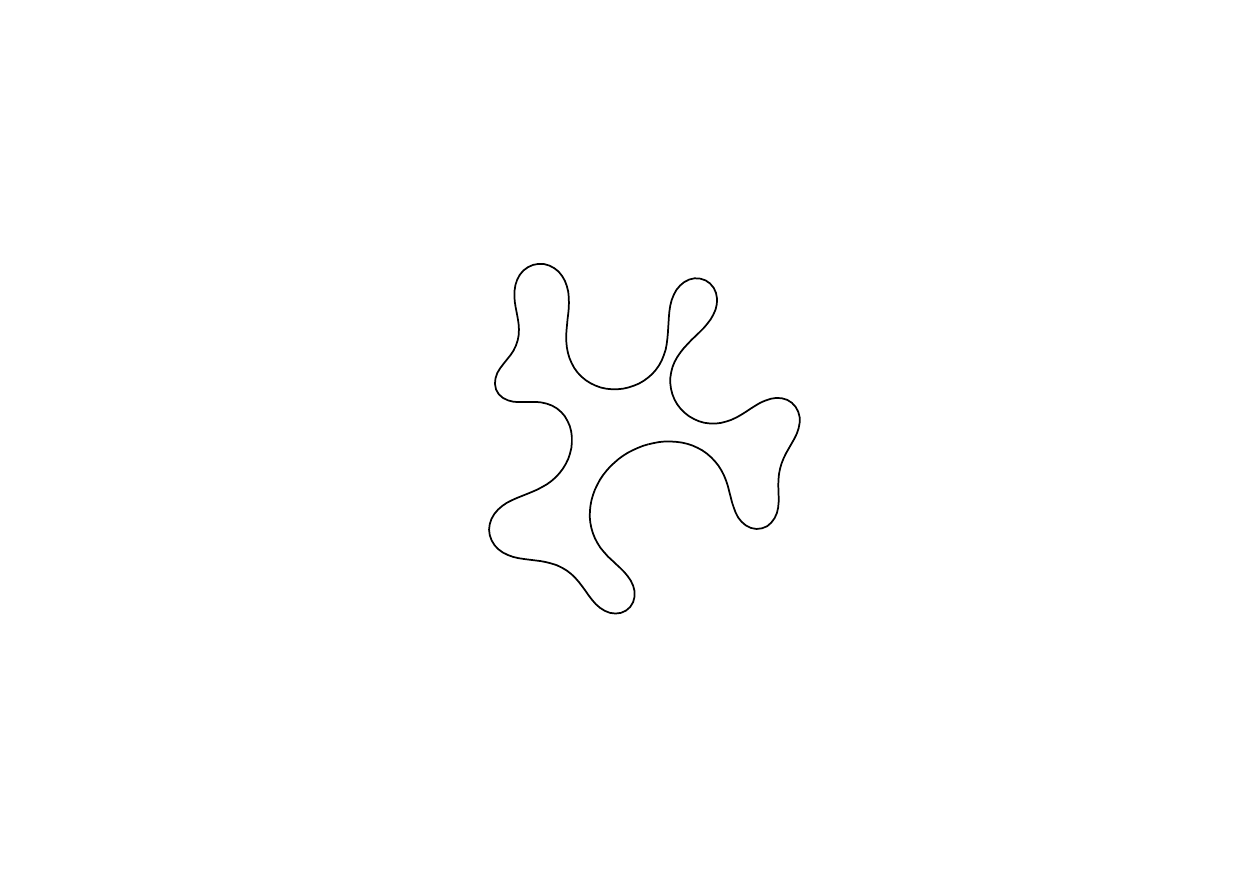}
                \subcaption{$t=1.84$}
            \end{minipage}
            \begin{minipage}{0.25\linewidth}
                \centering
                \includegraphics[scale = 0.32]{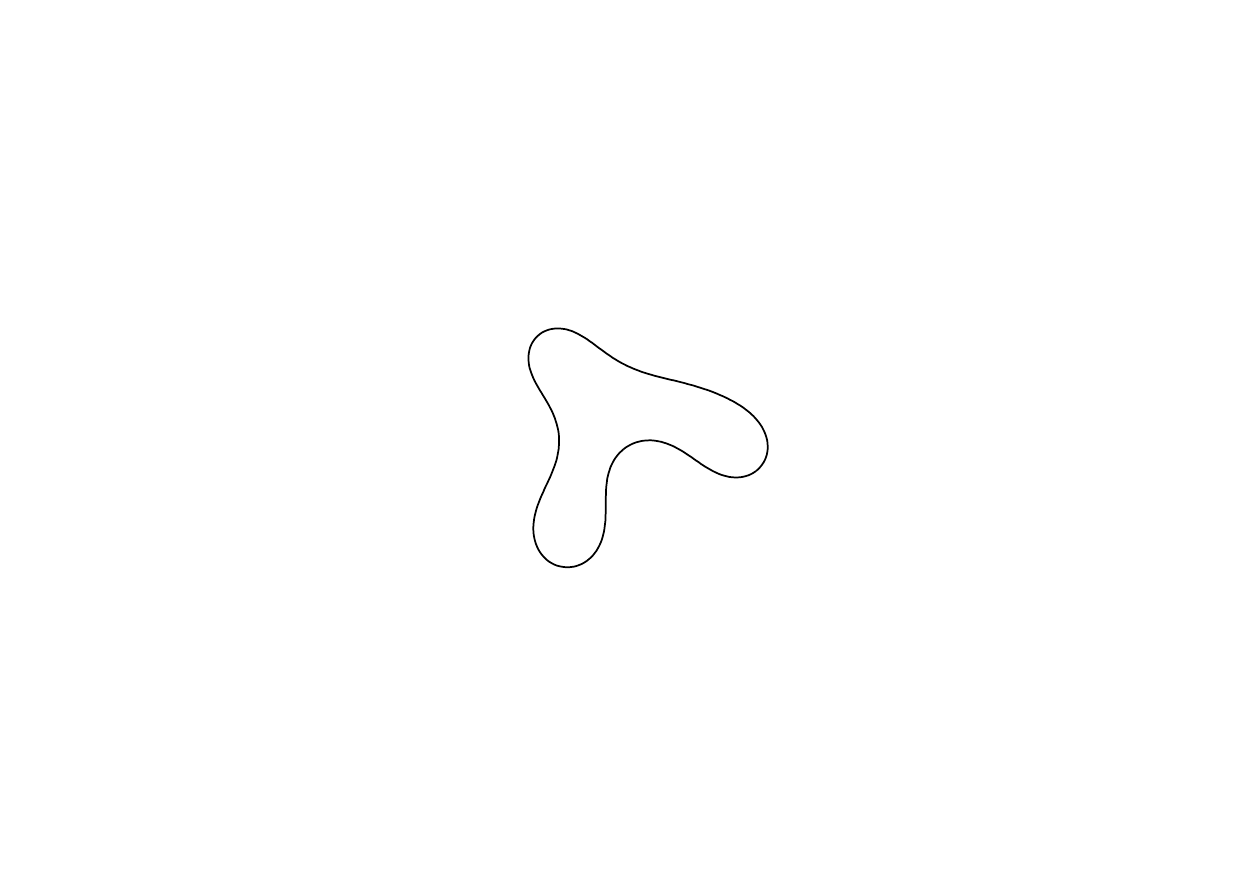}
                \subcaption{$t=2.21$}
            \end{minipage}
            \begin{minipage}{0.25\linewidth}
                \centering
                \includegraphics[scale = 0.32]{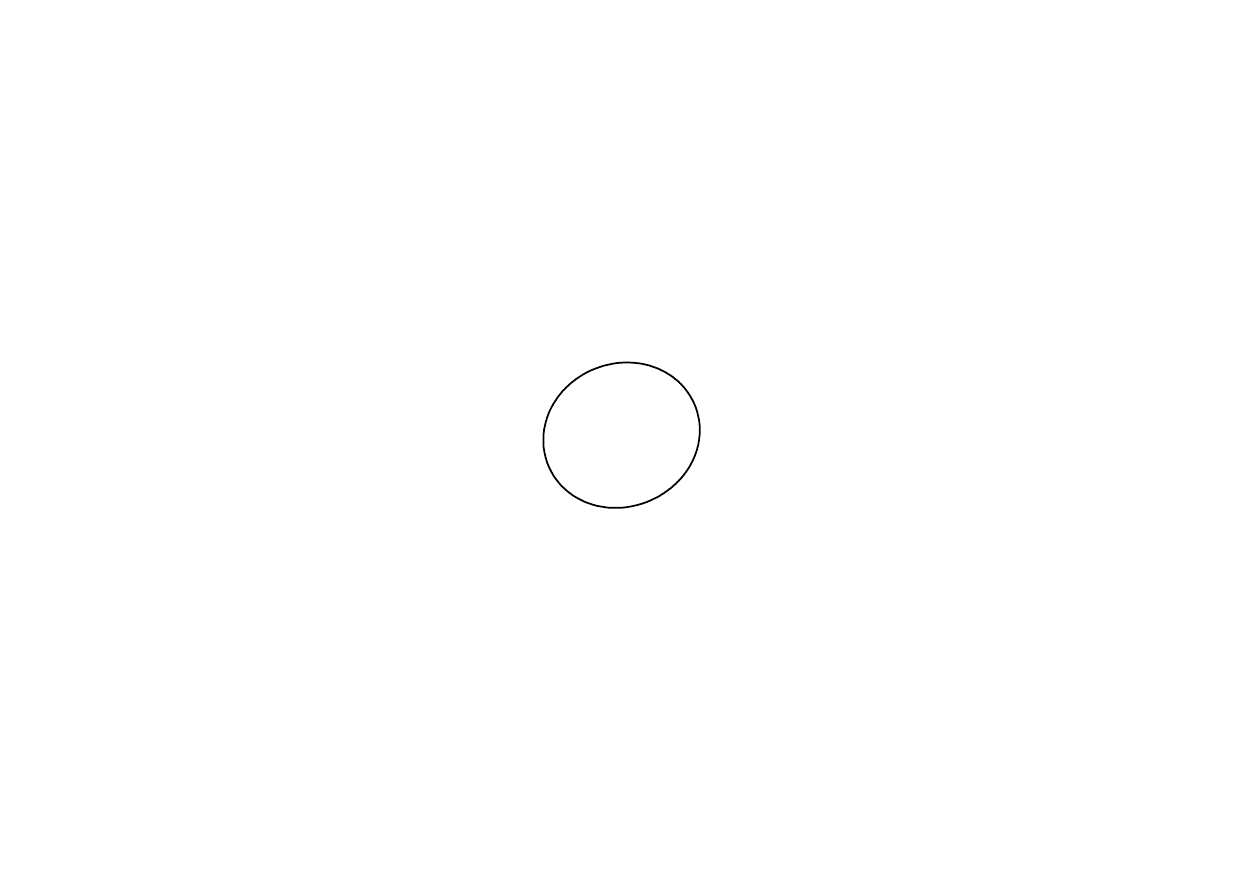}
                \subcaption{$t=2.58$}
            \end{minipage}
        \end{tabular}
        \caption{\label{fig:TDG}Results of numerical computation for Hele-Shaw flow with time-dependent gap: the time evolution of the boundary $\Gamma(t)$}
    \end{figure}

    The time-dependent gap $h(t)$ is given by $h(t) = \exp(t)$.
    In Fig.\ref{fig:TDG}, the time evolution of the boundary $\Gamma(t)$ is shown at $t = 0$, $0.36$, $0.73$, $1.10$, $1.47$, $1.84$, $2.21$, and $2.58$.
    The shapes at each time step almost coincide with the earlier study by Shelley et al.~\cite{shelley}.
    Note that this is not a comparison at the same discrete time because we use a different time discretization than in previous studies.
    
\begin{figure}[tb]
        \centering
        \includegraphics[scale = 0.5]{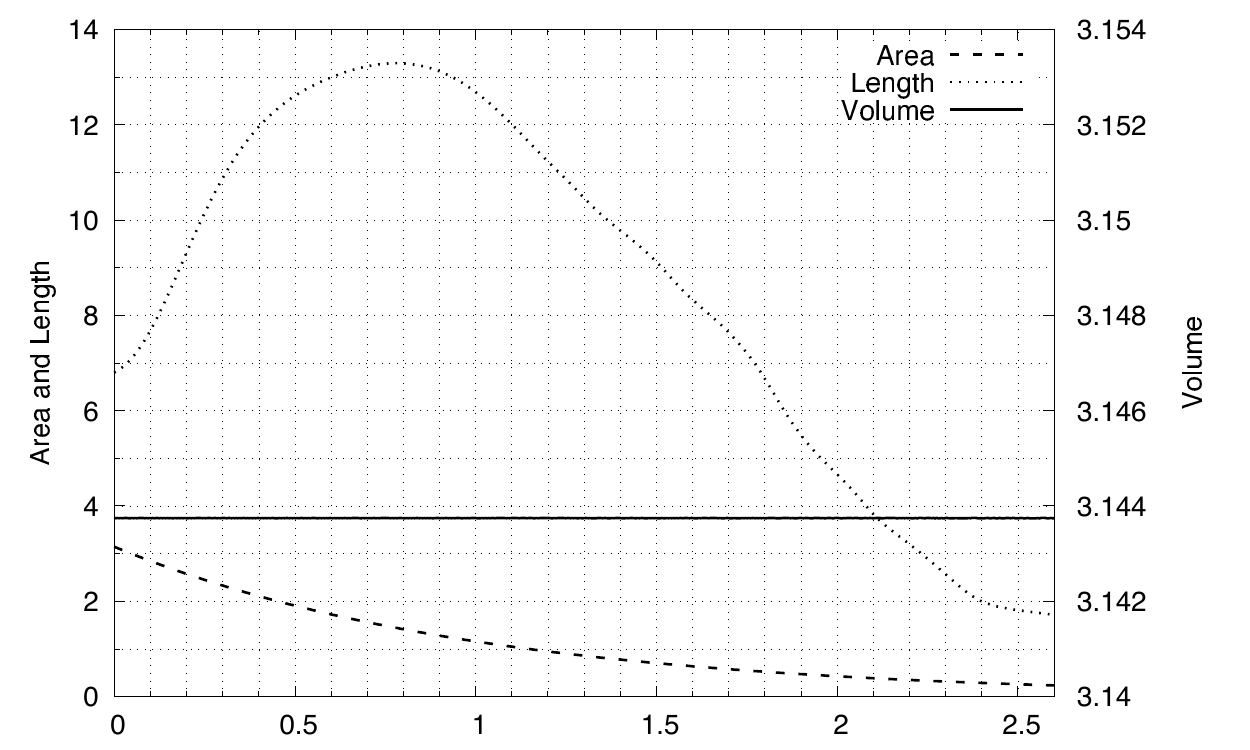}
        \caption{\label{fig:TDG_vol_preserving}Results of numerical computation for Hele-Shaw flow with time-dependent gap: time evolution of the total length, the area and the volume}    
    \end{figure}

    Fig.\ref{fig:TDG_vol_preserving} shows the time evolution of the total length, the area and the volume, where the horizontal axis represents the time, the left vertical axis represents the total length and the area, and the right vertical axis represents the volume.
    The graph shows that the area decreases monotonically and that the volume-preserving property is satisfied.
    
    
    \subsection{Hele-Shaw problem for magnetic fluids}
    
    In this subsection, we report the numerical results of the Hele-Shaw problem for magnetic fluids \eqref{eq:HS_magnetic}.
    
    \subsubsection{Computation of the magnetostatic field potential $\varphi_m$ by adaptive Monte Carlo method}
    
    To solve problem \eqref{eq:HS_magnetic}, we need to compute the magnetostatic field potential $\varphi_m$; that is, integration on $\Omega(t)$.
    In this paper, we employ the adaptive Monte Carlo method.
    
    At every time step, the region $\mathcal{MC}(t) \supseteq \Omega(t)$ is defined as follows:
    \begin{align*}
        \mathcal{MC}(t) 
        &\coloneqq 
        \left\{
            \begin{pmatrix}
                r(t)\cos\theta+g_x(t)\\
                r(t)\sin\theta+g_y(t)
            \end{pmatrix}
            \relmiddle|
            r(t)\in[0,r_{\max}(t)],\ \theta\in[0,2\pi)
        \right\},\\
        r_{\mathrm max}(t) &\coloneqq \max\{|\bm{X}_i(t) - \bm{g}(t)| \mid i = 1, 2, \ldots, N\},
    \end{align*}
    where $\bm{g}(t) = (g_x(t), g_y(t))^{\top}$ is the barycenter of $\Omega(t)$.
    Then, $M$ sample points $\{\bm{p}_k\}_{k=1}^M$ are distributed in $\mathcal{MC}$ uniformly.
    
    Distributed sample points $\{\bm{p}_k\}_{k=1}^M$ are judged whether $\bm{p}_k \in \Omega$ by the following algorithm.
    Let $a_k^i$ be the angle between the vectors $\bm{X}_{i-1} - \bm{p}_k$ and $\bm{X}_i - \bm{p}_k$, i.e.,
    \begin{equation*}
        a_k^i=\arccos\left(\frac{\bm{X}_{i-1}-\bm{p}_k}{|\bm{X}_{i-1}-\bm{p}_k|}\cdot\frac{\bm{X}_i-\bm{p}_k}{|\bm{X}_i-\bm{p}_k|}\right).
    \end{equation*}
    If $\sum_{i=1}^Na_k^i=2\pi$ holds, then the sample point $\bm{p}_k$ is judged to be inside $\Omega$; otherwise, it is judged outside $\Omega$.
    
    Let $M_{\mathrm in}$ be the number of points $\bm{p}_k$ belonging to $\Omega$. Then the area element $dS$ of $\Omega$ is defined as $dS = A / M_{in}$. By using this $dS$, the the value of $\varphi = \varphi_m / M_c$ at the vertex $\bm{X}_i$ can be computed approximately as follows:
    \begin{equation*}
        \varphi(\bm{X}_i) = -\sum_{\bm{p}_k \in \Omega} \left(\dfrac{1}{|\bm{X}_i - \bm{p}_k|} - \dfrac{1}{\sqrt{|\bm{X}_i - \bm{p}_k|^2 + h^2}}\right)dS
    \end{equation*}
    for $i = 1, 2, \ldots, N$.
    
    \subsubsection{Numerical results}
    
\begin{figure}[tb]
        \begin{tabular}{cccc}
            \begin{minipage}{0.25\linewidth}
                \centering
                \includegraphics[scale = 0.32]{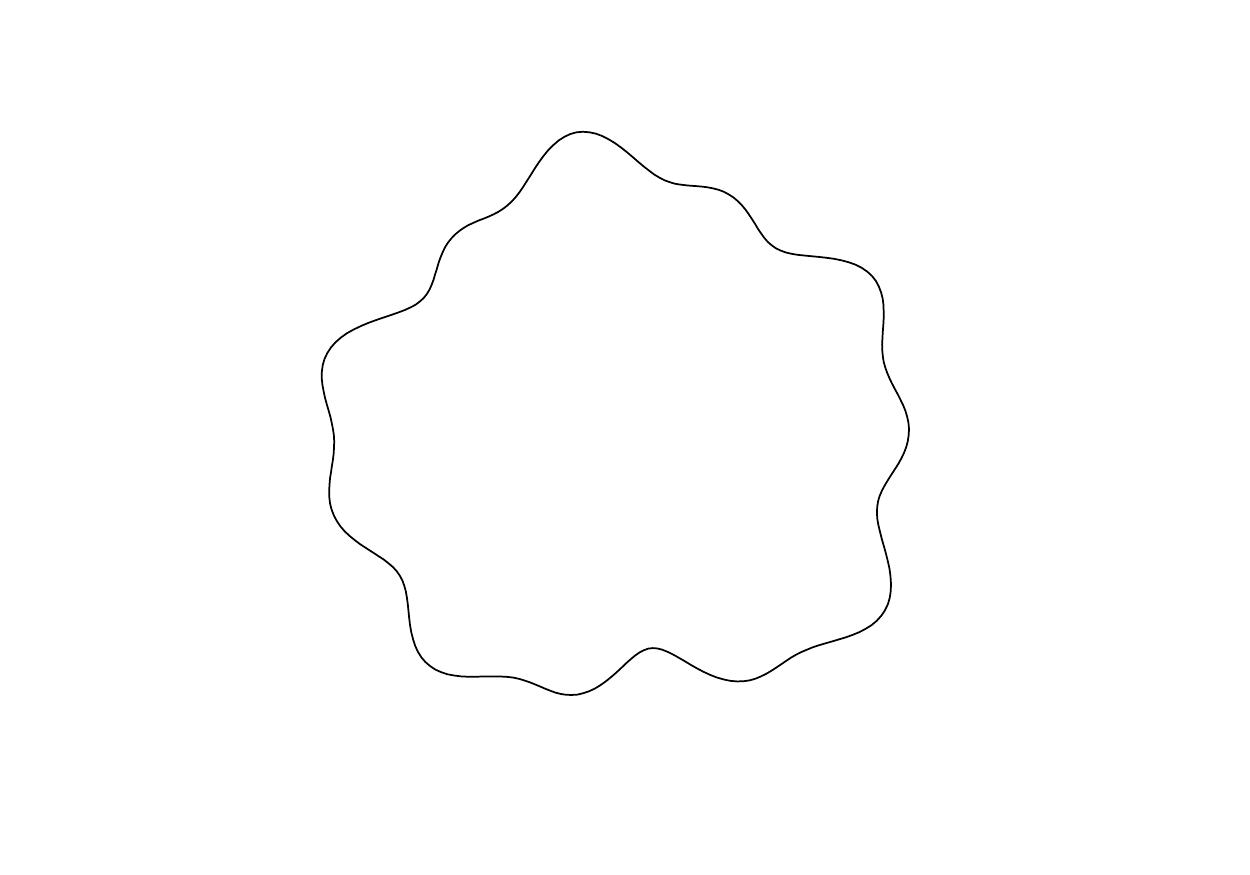}
                \subcaption{$t=0.0$}
            \end{minipage}
            \begin{minipage}{0.25\linewidth}
                \centering
                \includegraphics[scale = 0.32]{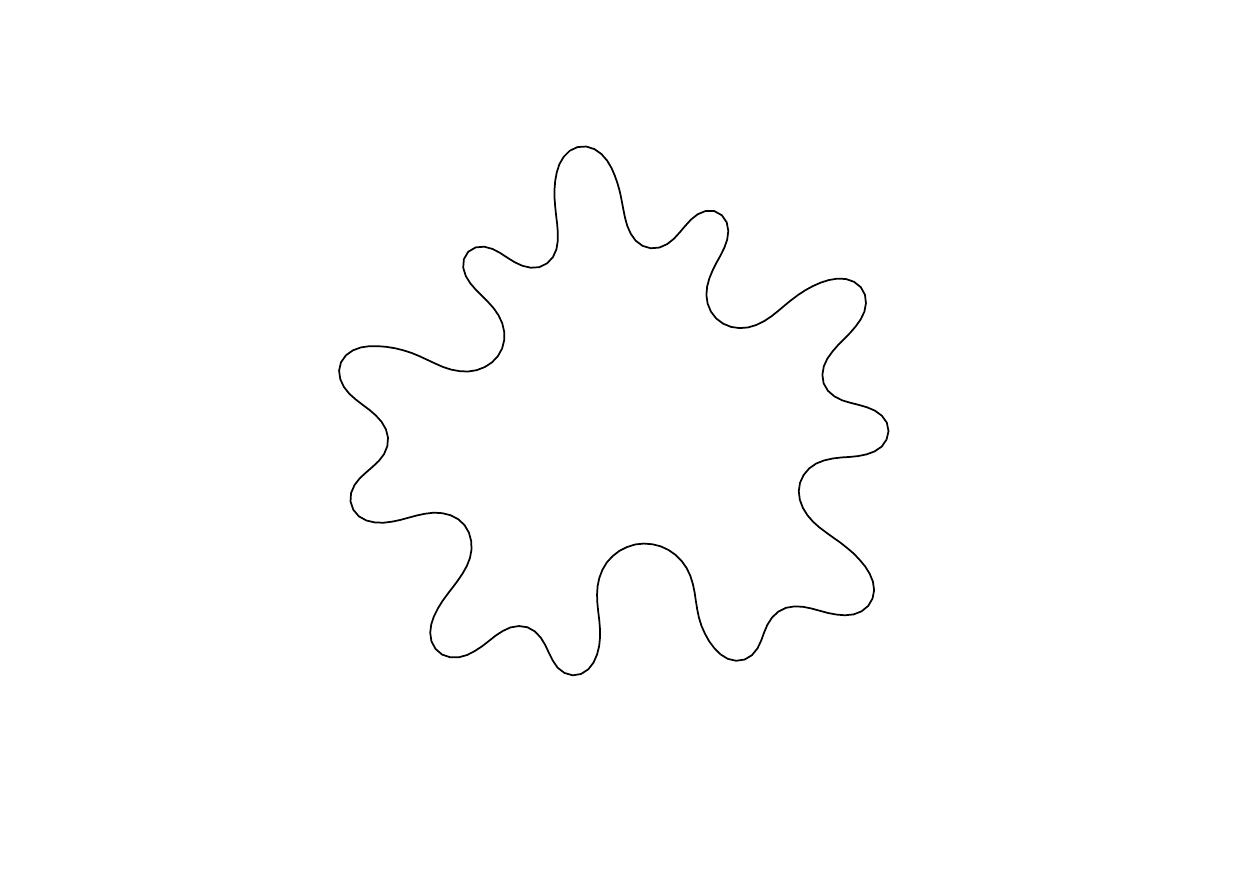}
                \subcaption{$t=0.42$}
            \end{minipage}
            \begin{minipage}{0.25\linewidth}
                \centering
                \includegraphics[scale = 0.32]{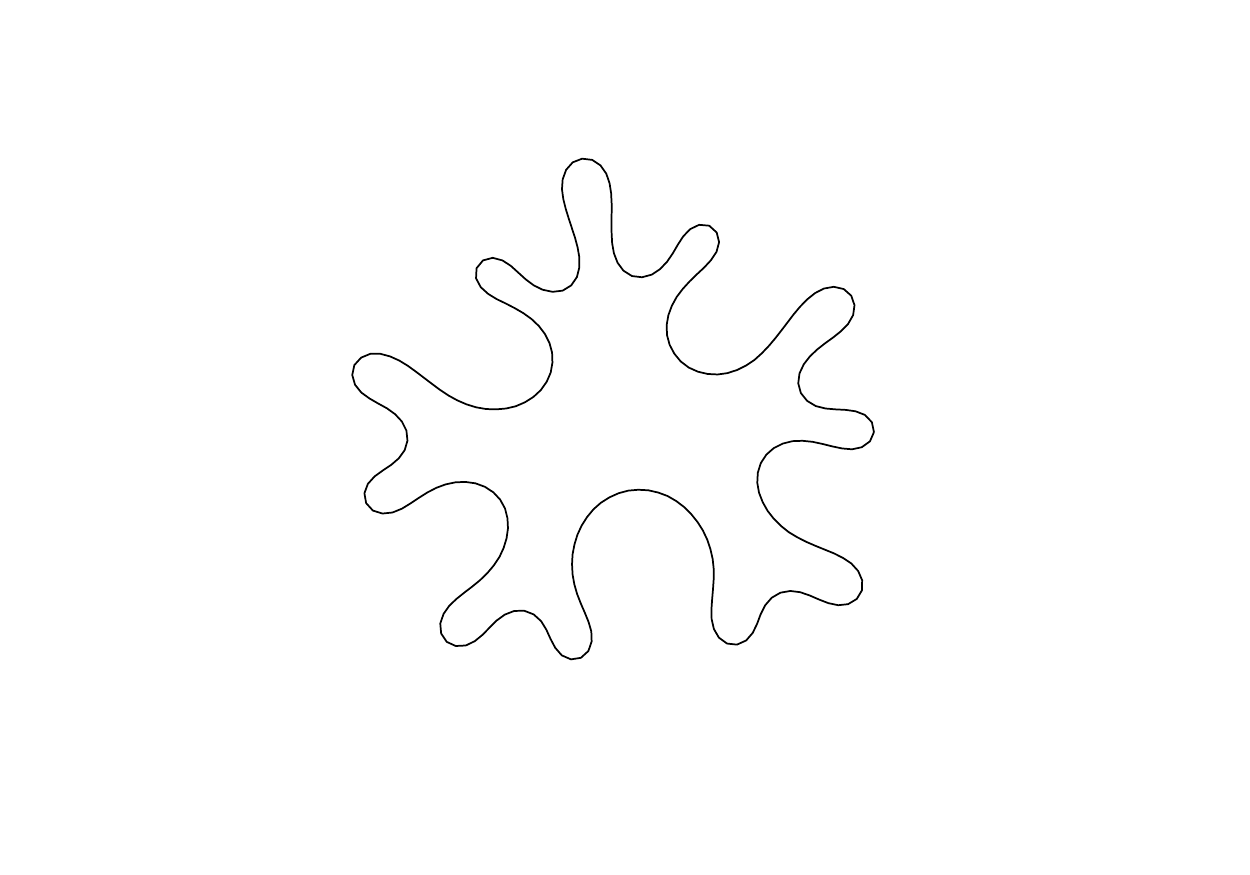}
                \subcaption{$t=0.84$}
            \end{minipage}
            \begin{minipage}{0.25\linewidth}
                \centering
                \includegraphics[scale = 0.32]{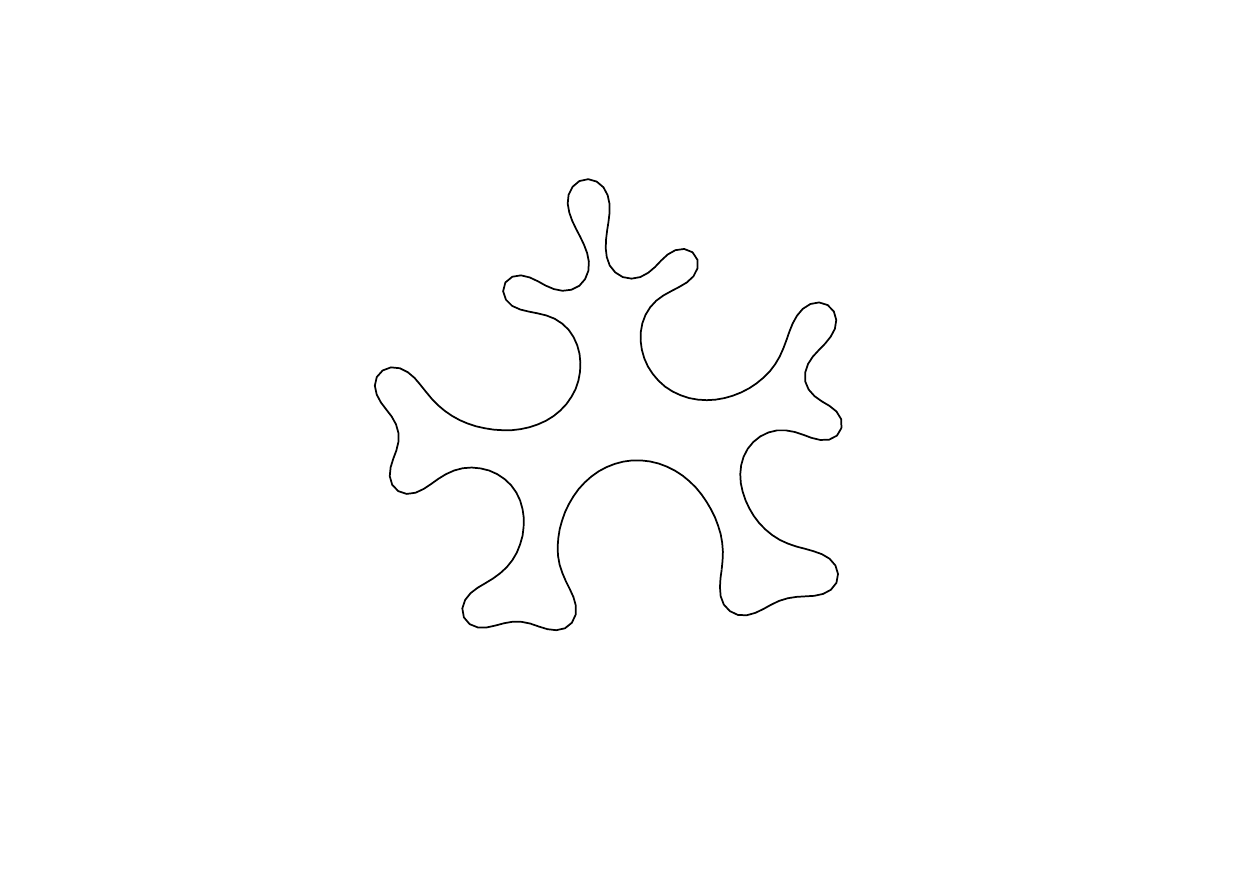}
                \subcaption{$t=1.26$}
            \end{minipage}
            \\
            \begin{minipage}{0.25\linewidth}
                \centering
                \includegraphics[scale = 0.32]{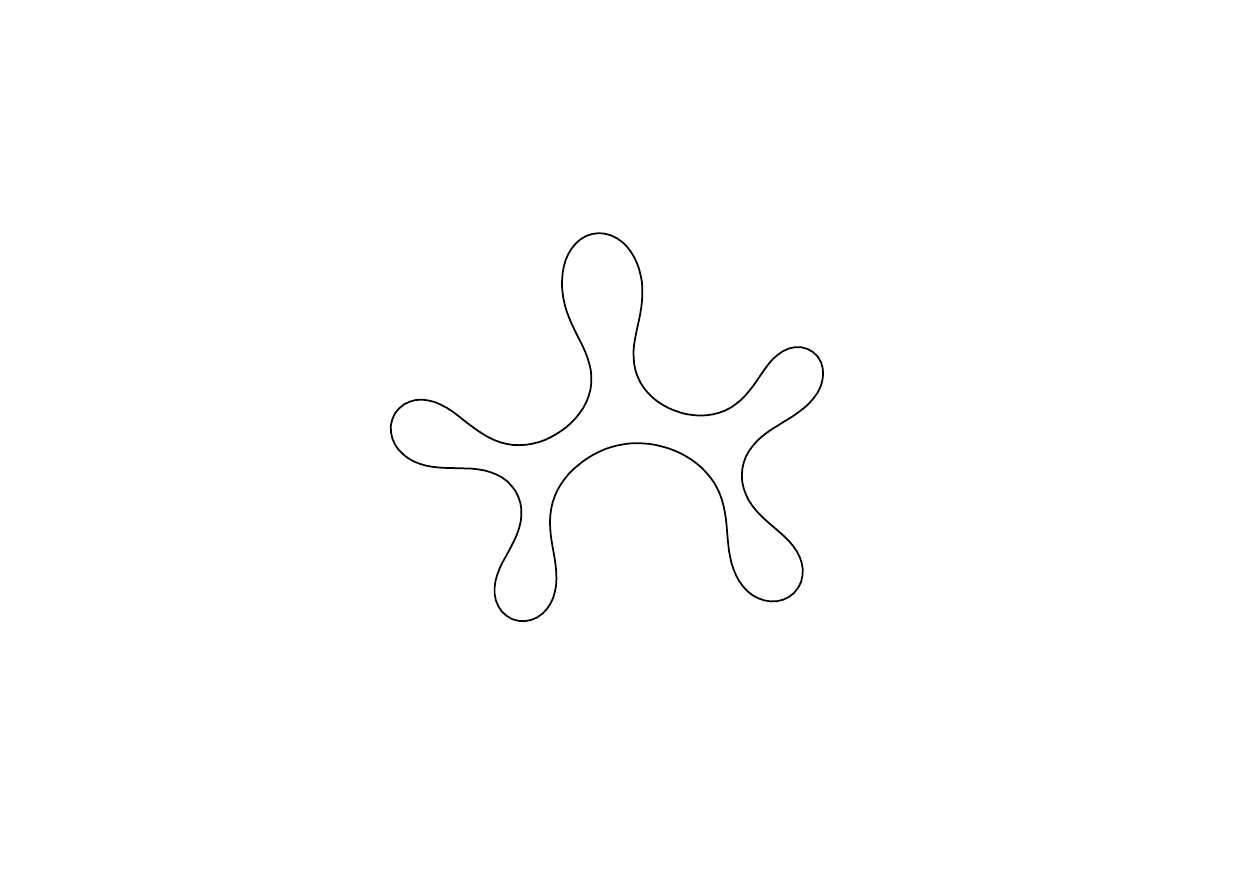}
                \subcaption{$t=1.68$}
            \end{minipage}
            \begin{minipage}{0.25\linewidth}
                \centering
                \includegraphics[scale = 0.32]{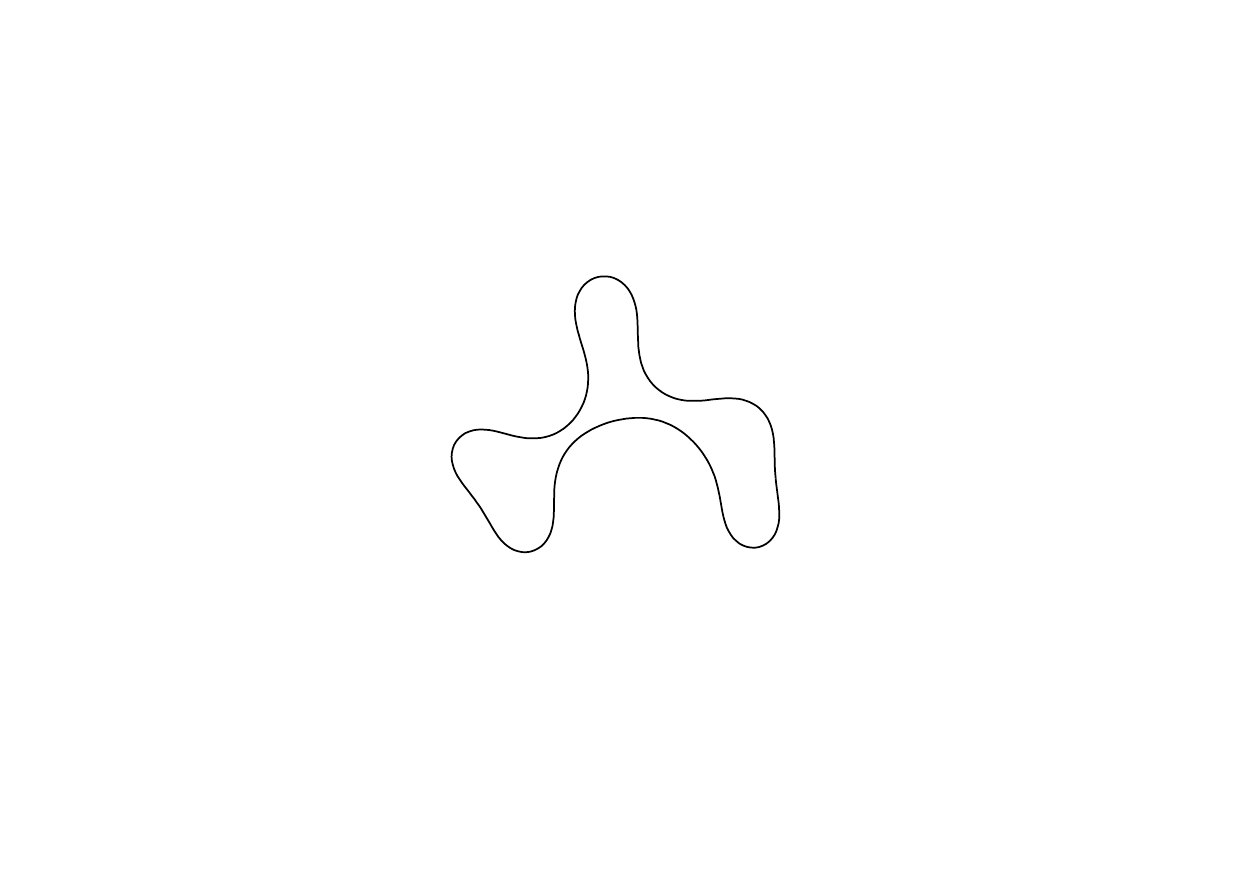}
                \subcaption{$t=2.10$}
            \end{minipage}
            \begin{minipage}{0.25\linewidth}
                \centering
                \includegraphics[scale = 0.32]{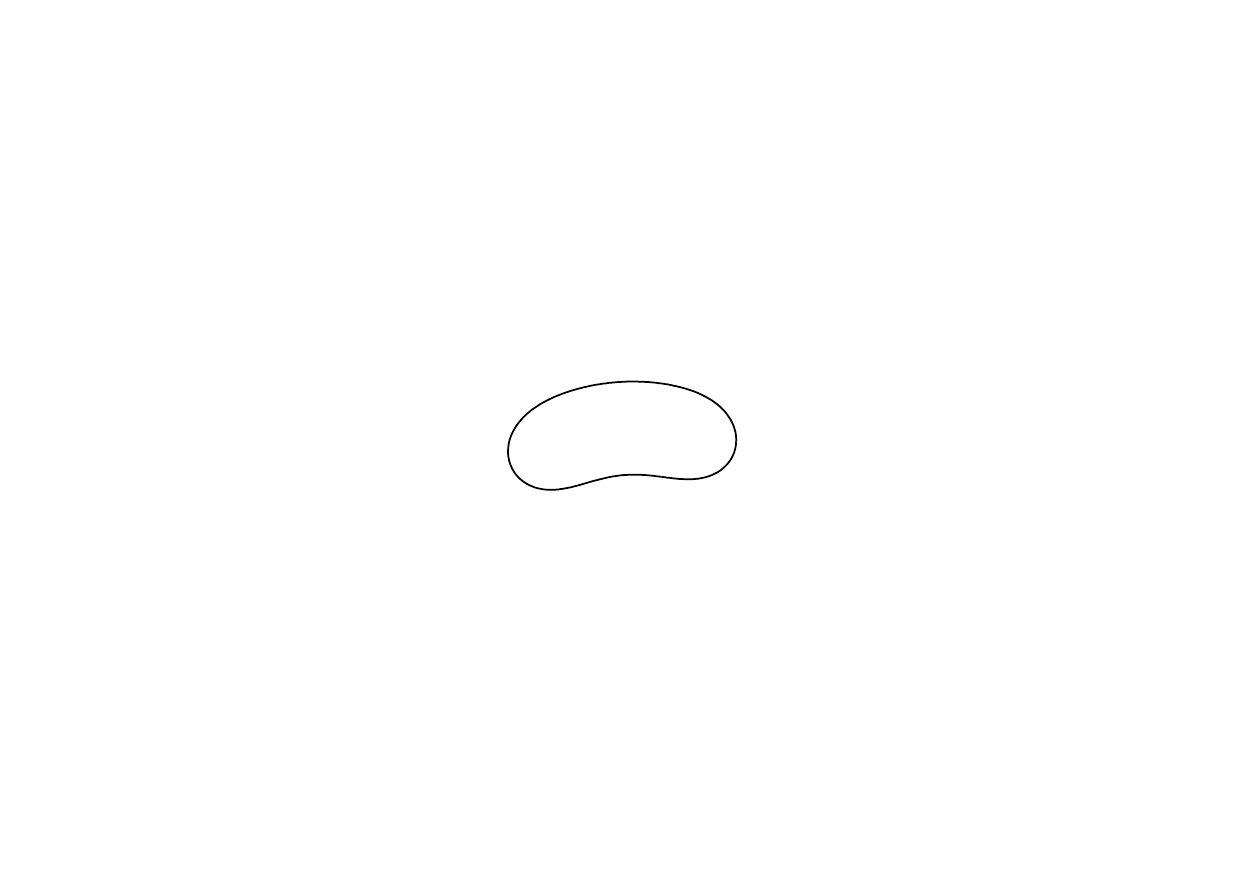}
                \subcaption{$t=2.52$}
            \end{minipage}
            \begin{minipage}{0.25\linewidth}
                \centering
                \includegraphics[scale = 0.32]{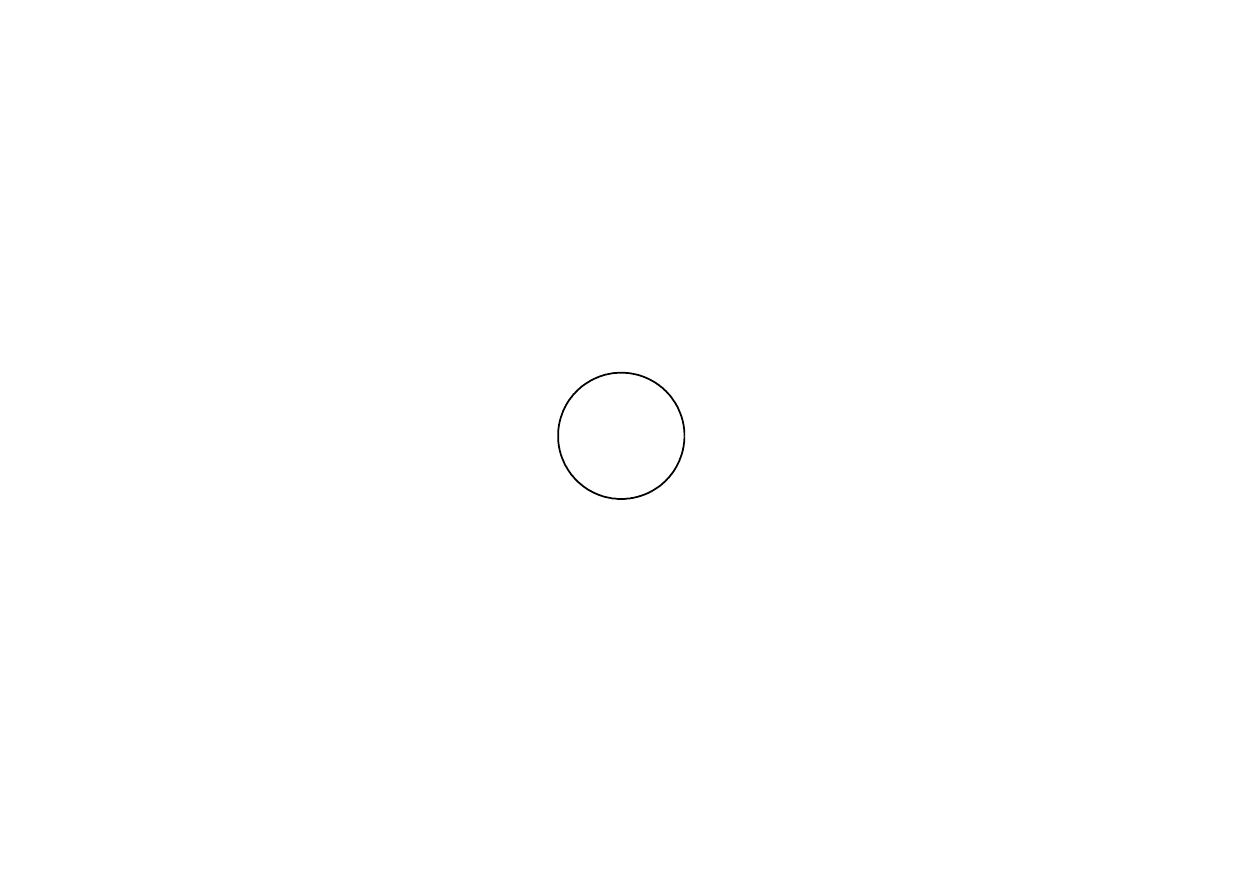}
                \subcaption{$t=2.94$}
            \end{minipage}
        \end{tabular}
            \caption{\label{fig:Bmv0Ca100}$\mathrm{Bmv} = 0, \mathrm{Ca} = 100$}
    \end{figure}

\begin{figure}[tb]
        \begin{tabular}{cccc}
            \begin{minipage}{0.25\linewidth}
                \centering
                \includegraphics[scale = 0.32]{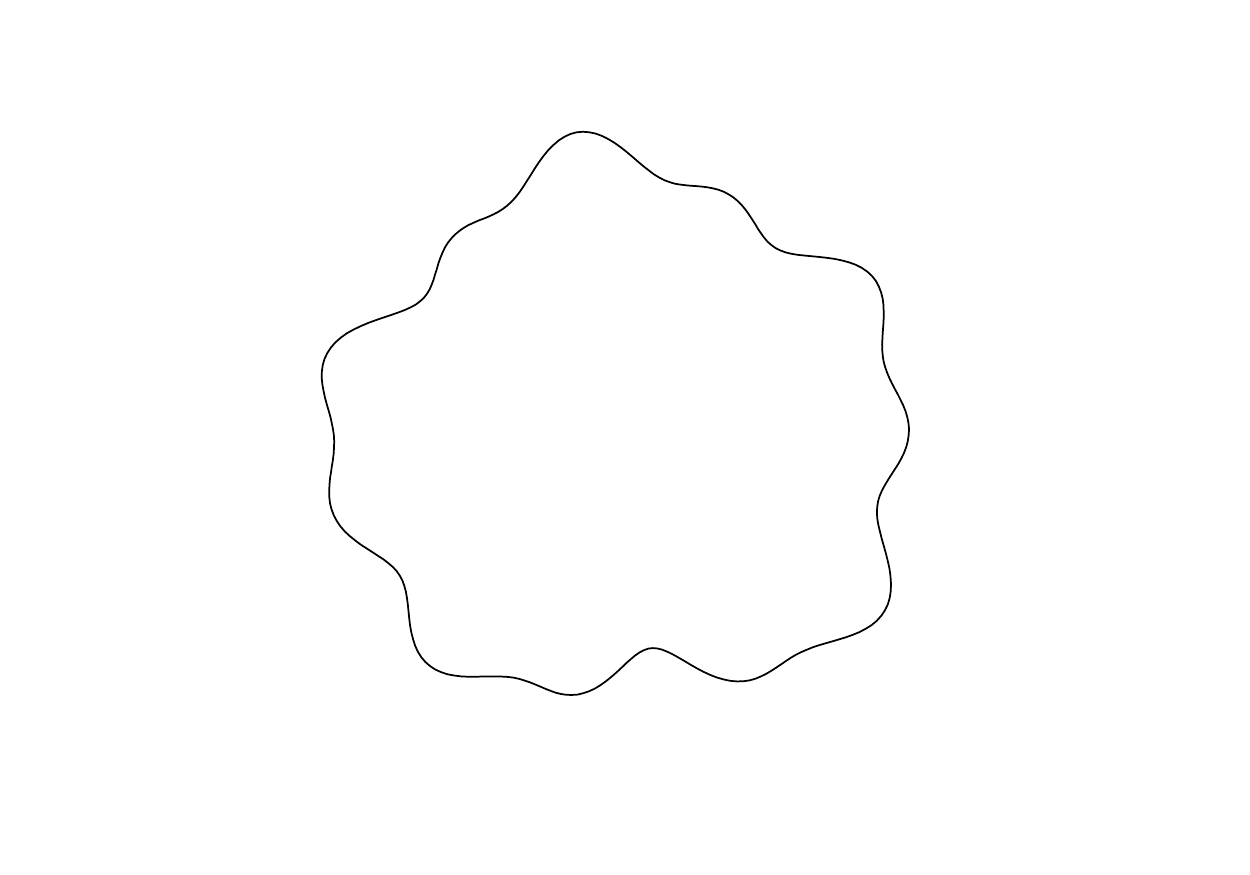}
                \subcaption{$t=0.0$}
            \end{minipage}
            \begin{minipage}{0.25\linewidth}
                \centering
                \includegraphics[scale = 0.32]{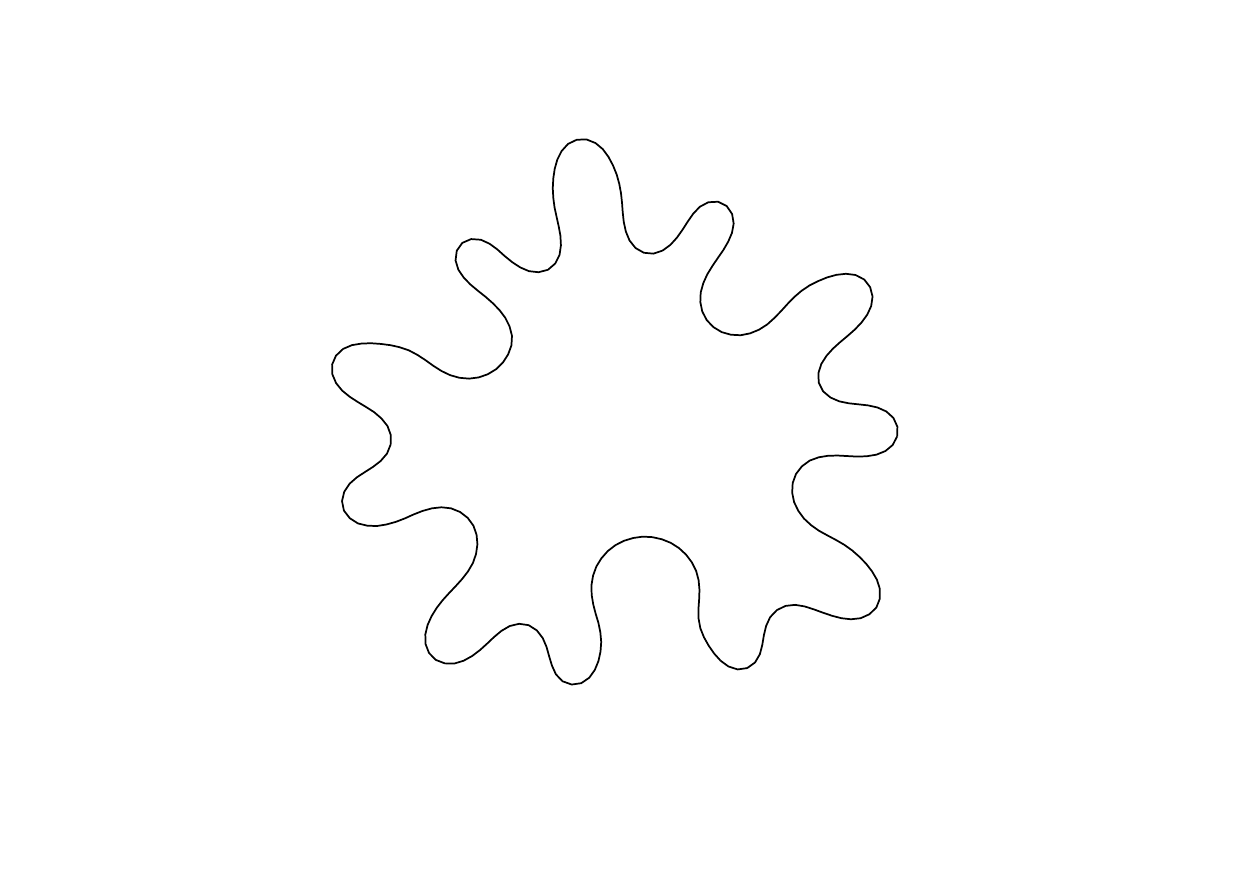}
                \subcaption{$t=0.42$}
            \end{minipage}
            \begin{minipage}{0.25\linewidth}
                \centering
                \includegraphics[scale = 0.32]{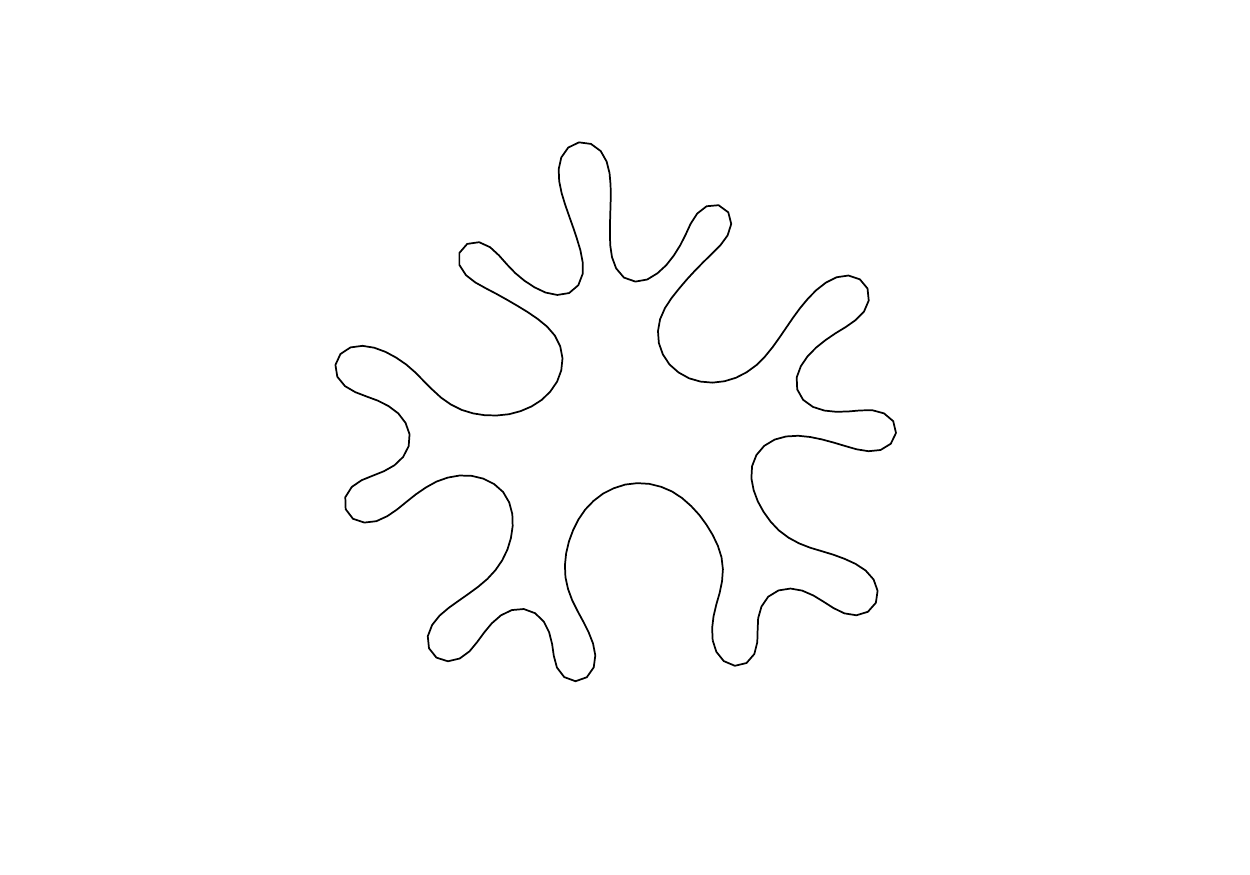}
                \subcaption{$t=0.84$}
            \end{minipage}
            \begin{minipage}{0.25\linewidth}
                \centering
                \includegraphics[scale = 0.32]{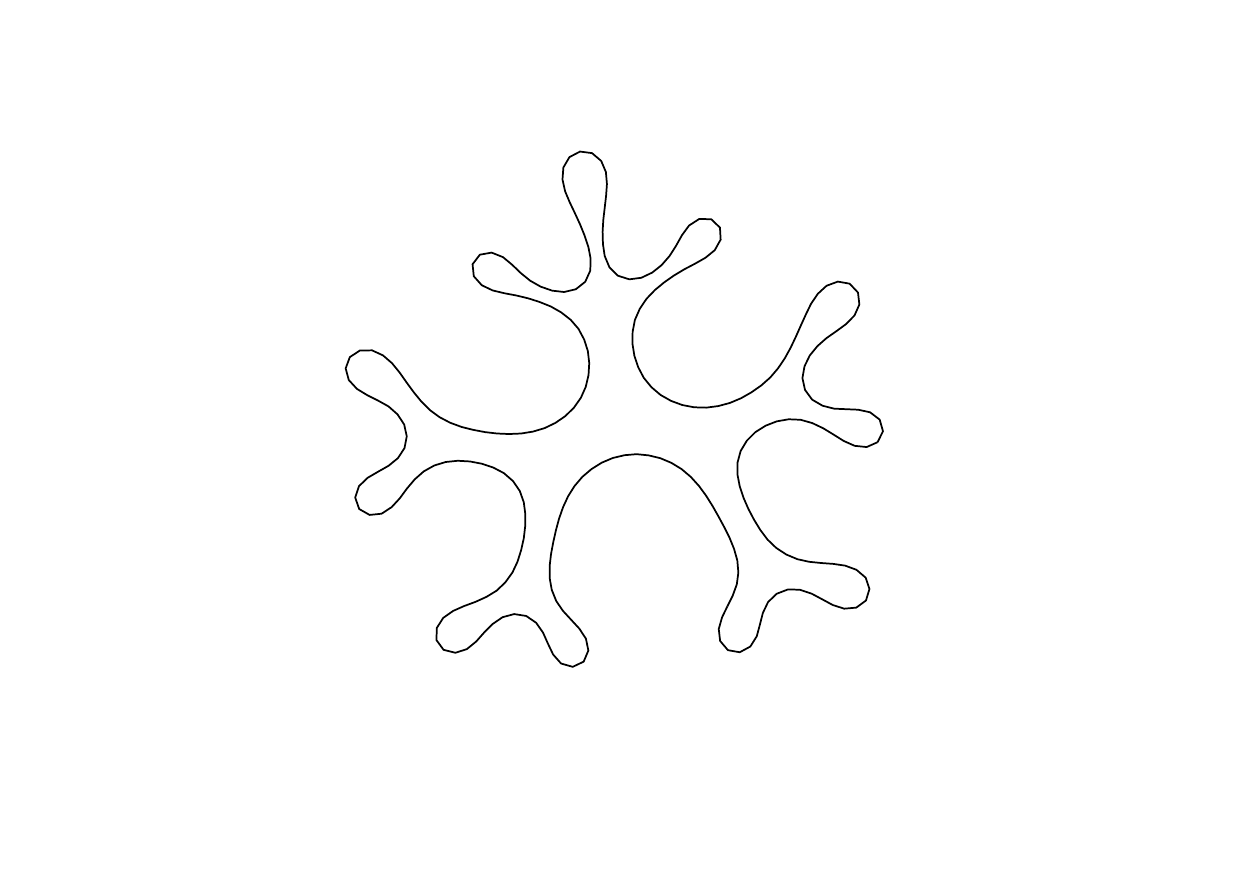}
                \subcaption{$t=1.26$}
            \end{minipage}
            \\
            \begin{minipage}{0.25\linewidth}
                \centering
                \includegraphics[scale = 0.32]{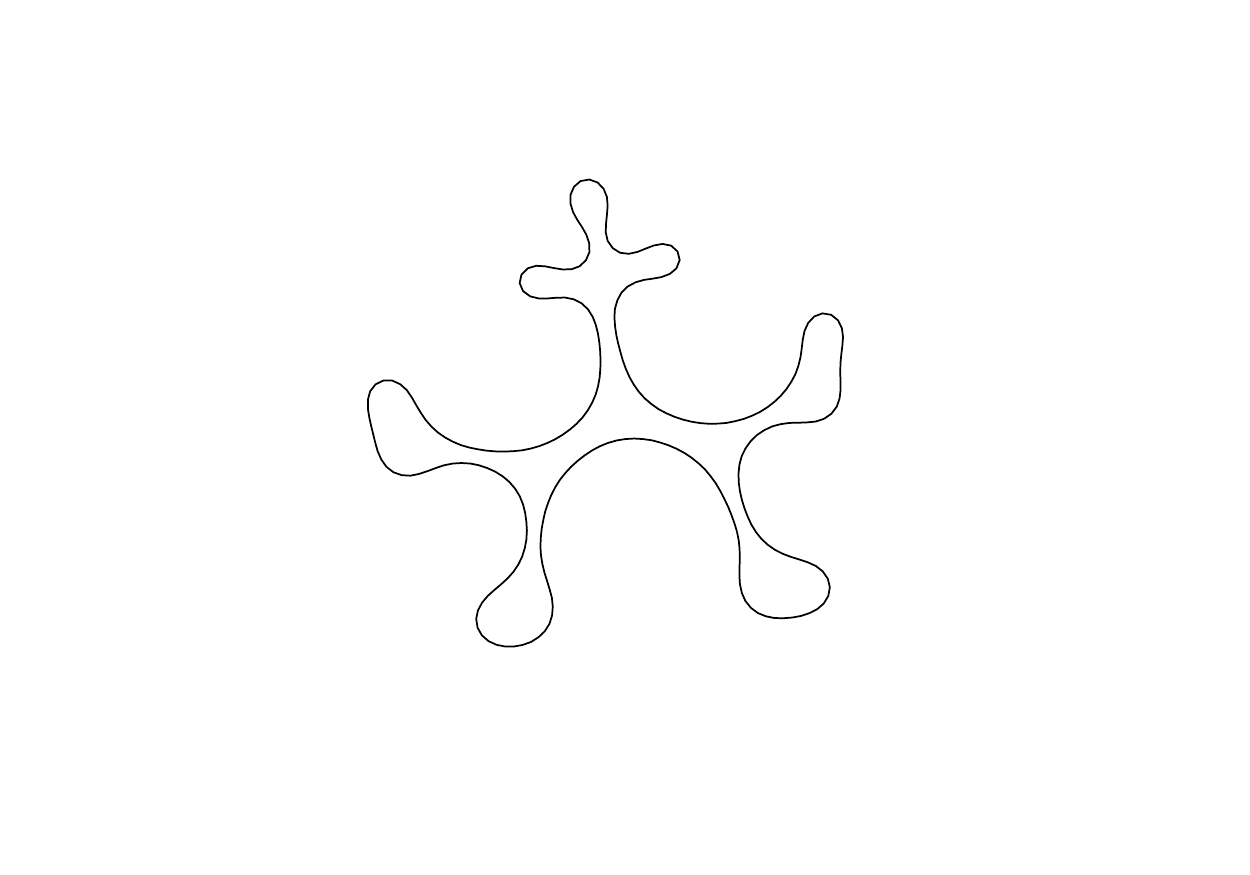}
                \subcaption{$t=1.68$}
            \end{minipage}
            \begin{minipage}{0.25\linewidth}
                \centering
                \includegraphics[scale = 0.32]{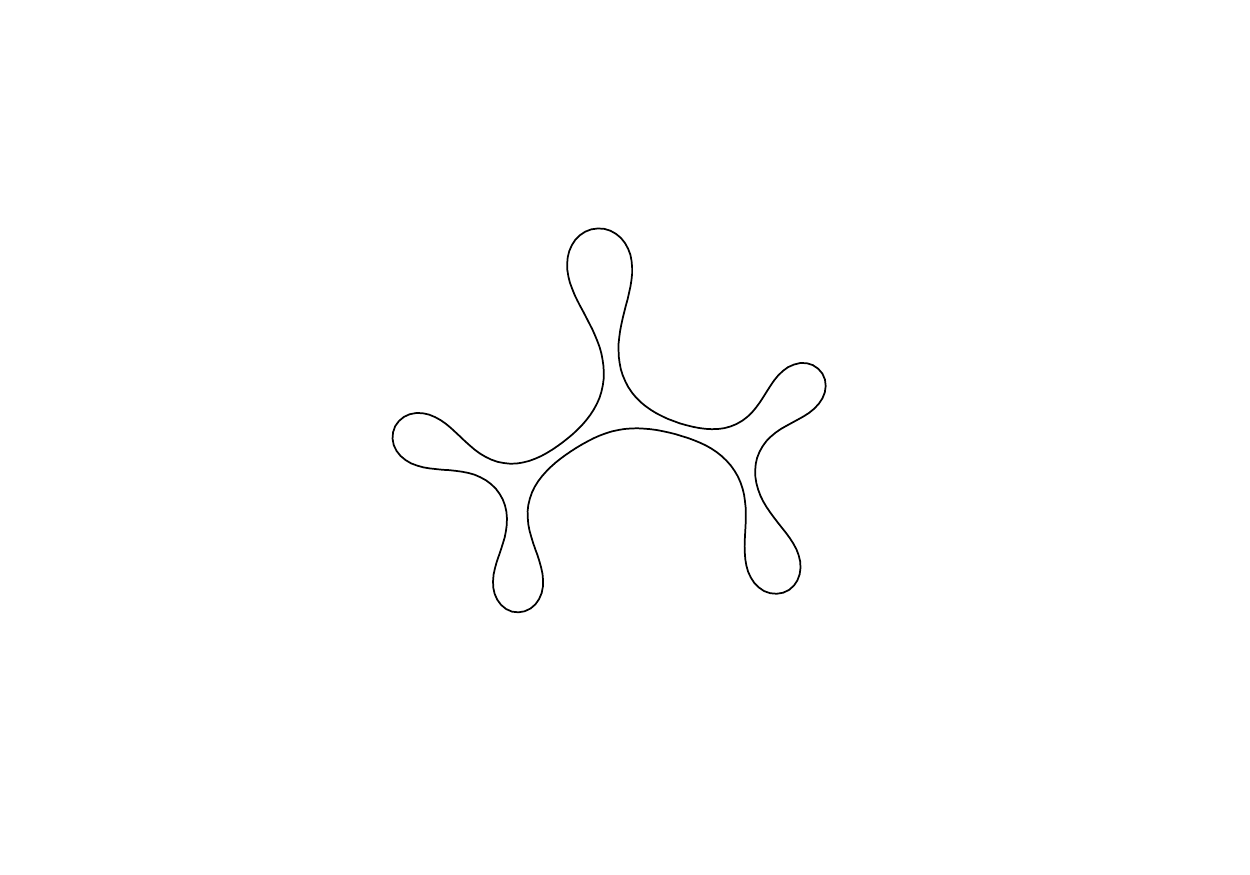}
                \subcaption{$t=2.10$}
            \end{minipage}
            \begin{minipage}{0.25\linewidth}
                \centering
                \includegraphics[scale = 0.32]{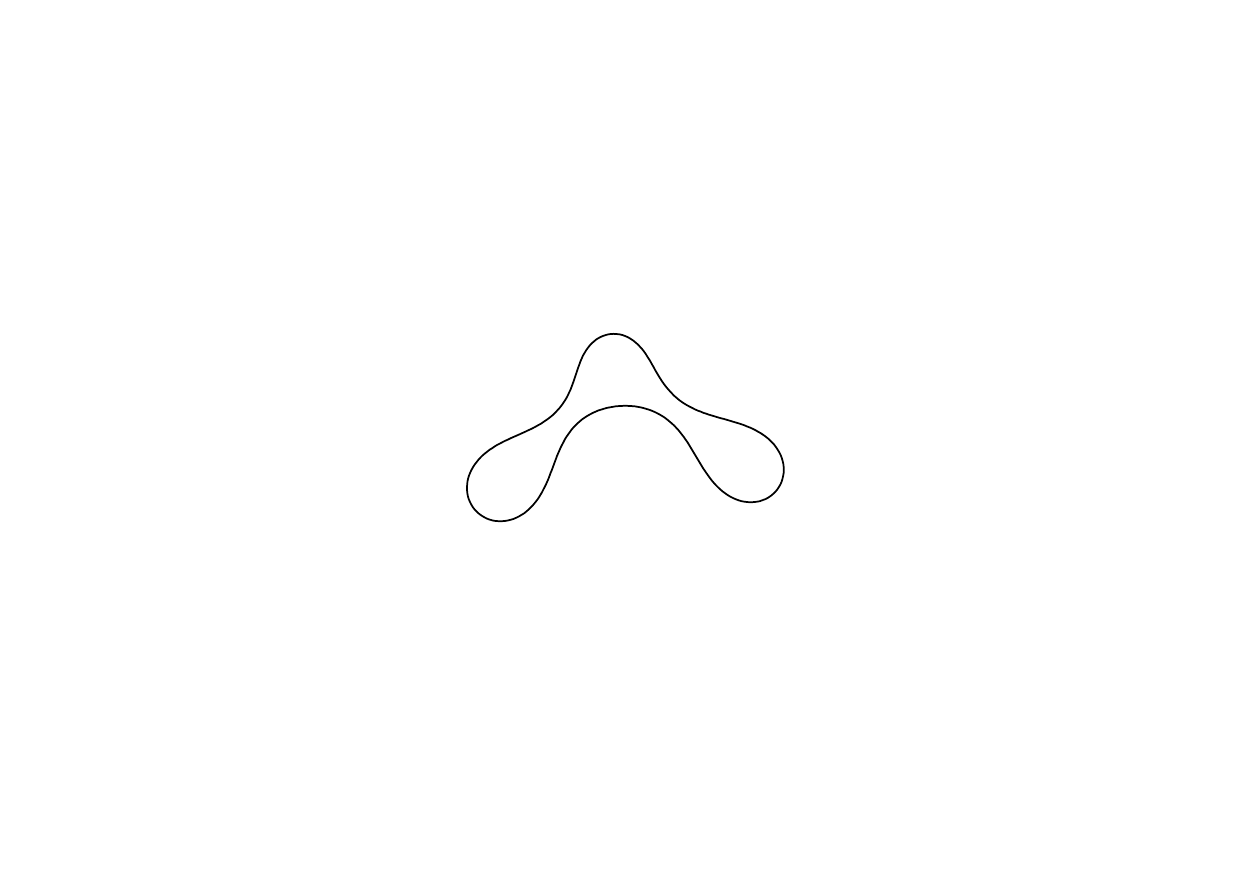}
                \subcaption{$t=2.52$}
            \end{minipage}
            \begin{minipage}{0.25\linewidth}
                \centering
                \includegraphics[scale = 0.32]{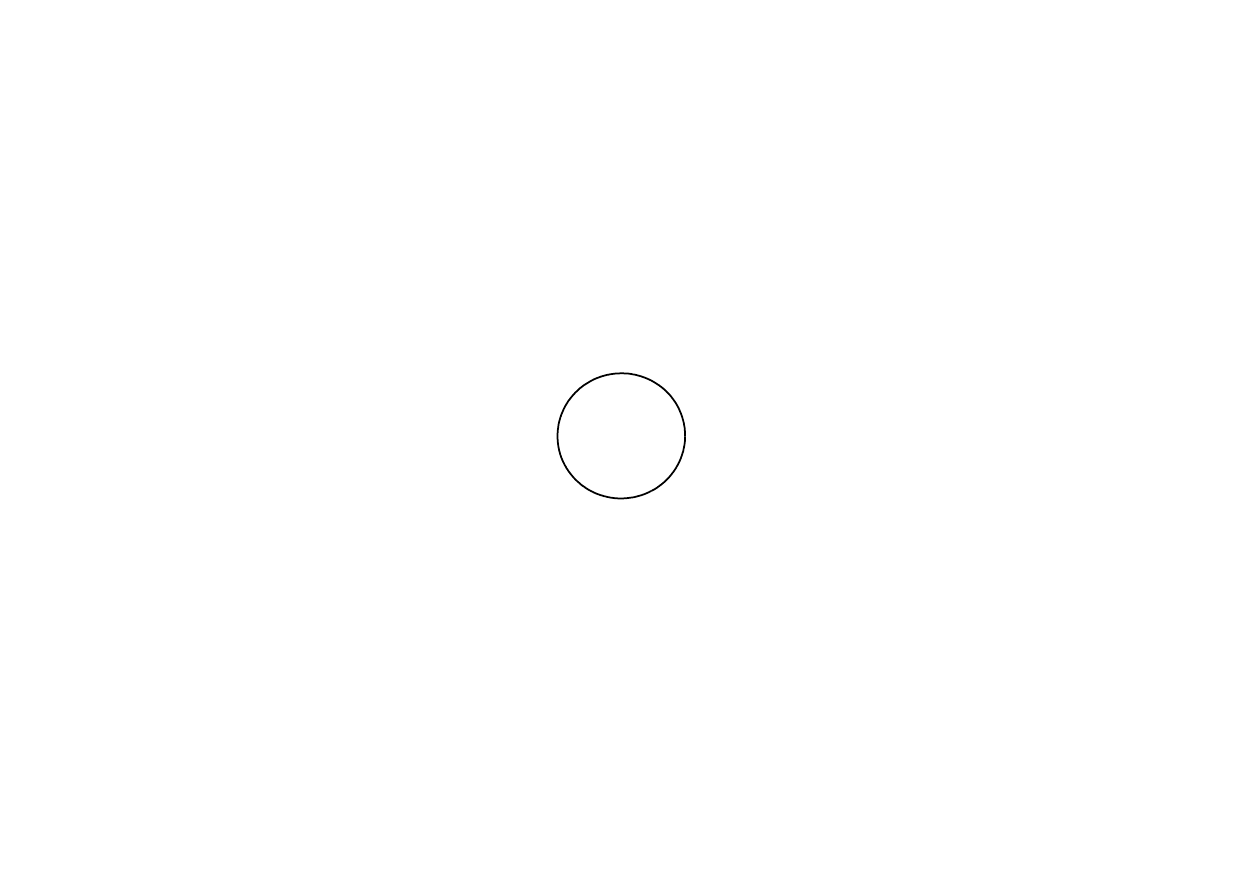}
                \subcaption{$t=2.94$}
            \end{minipage}
        \end{tabular}
        \caption{\label{fig:Bmv25Ca100}$\mathrm{Bmv} = 25, \mathrm{Ca} = 100$}
    \end{figure}

\begin{figure}[tb]
        \begin{tabular}{cccc}
            \begin{minipage}{0.25\linewidth}
                \centering
                \includegraphics[scale = 0.32]{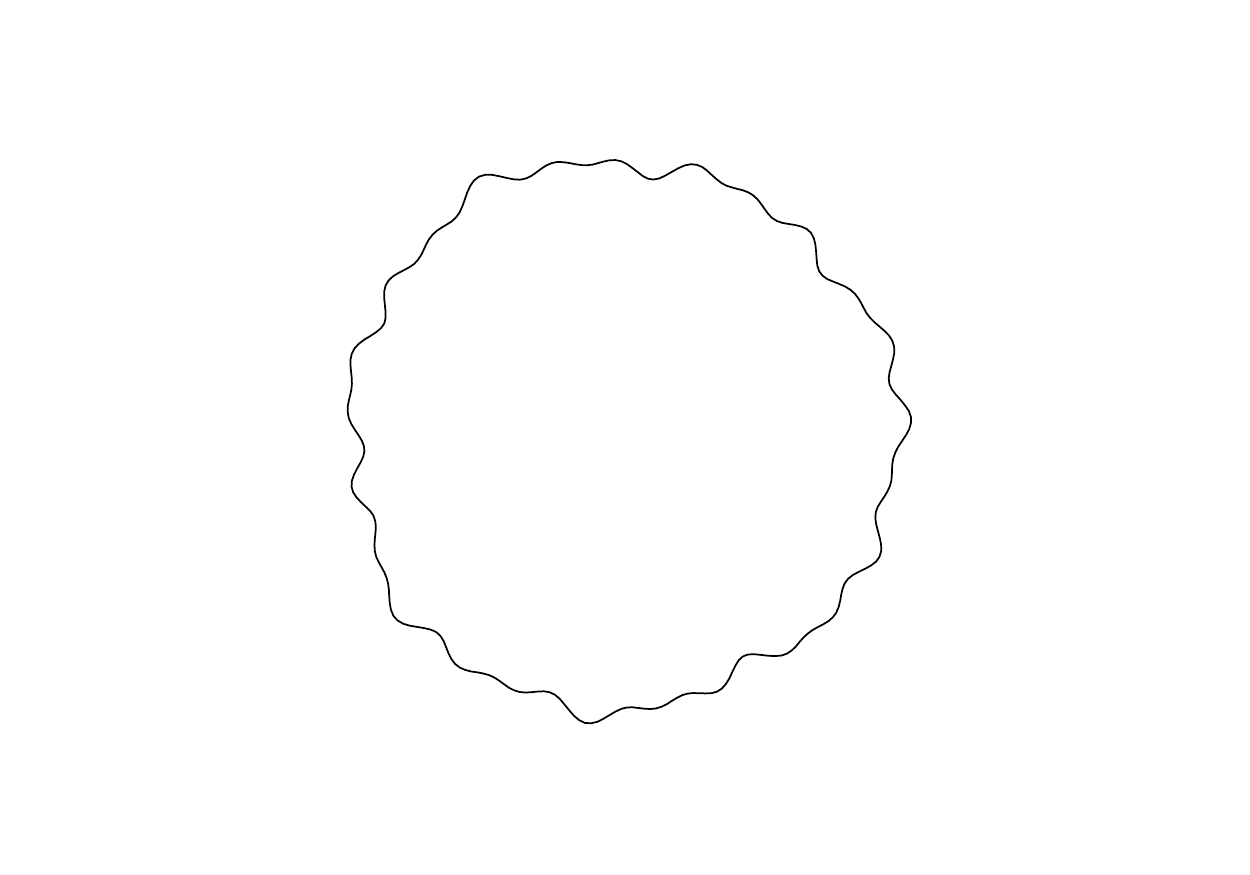}
                \subcaption{$t=0.0$}
            \end{minipage}
            \begin{minipage}{0.25\linewidth}
                \centering
                \includegraphics[scale = 0.32]{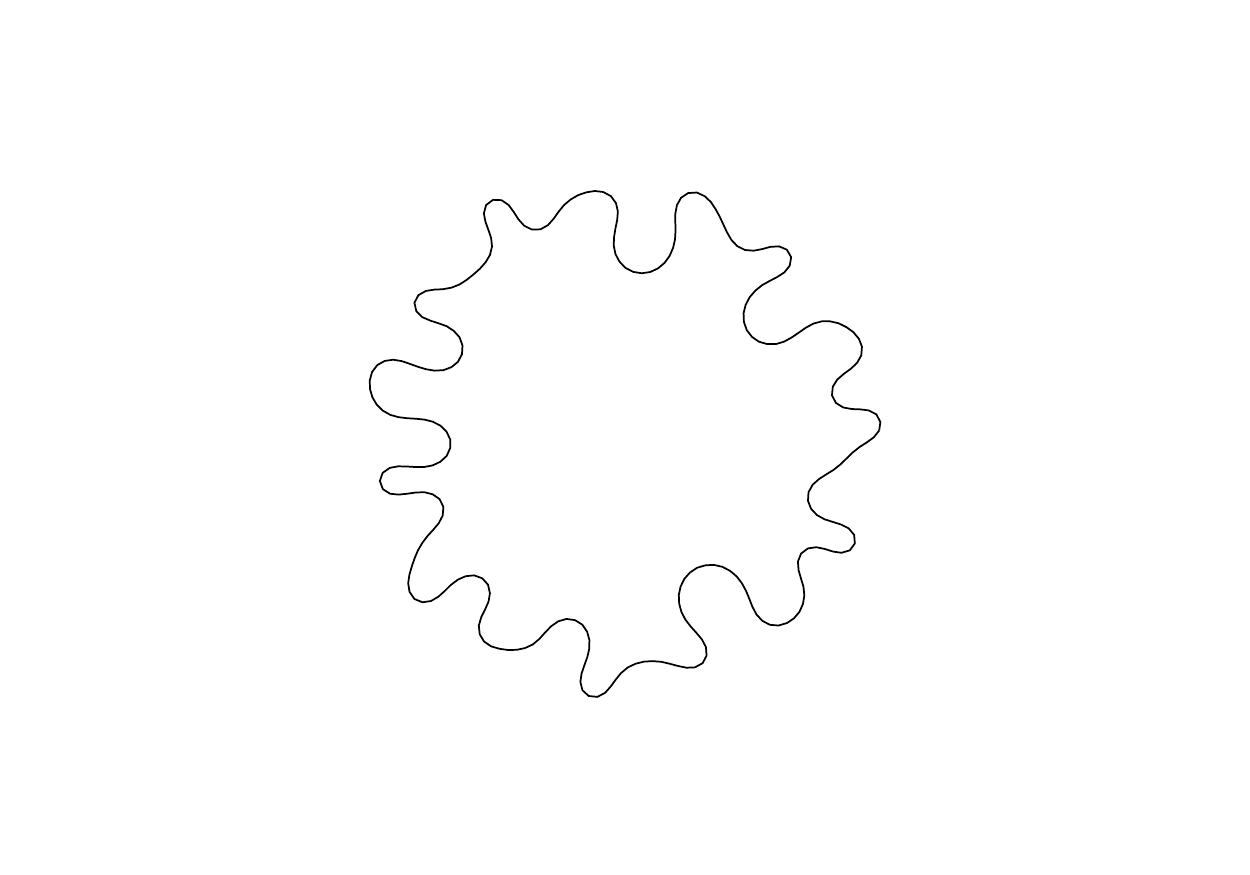}
                \subcaption{$t=0.42$}
            \end{minipage}
            \begin{minipage}{0.25\linewidth}
                \centering
                \includegraphics[scale = 0.32]{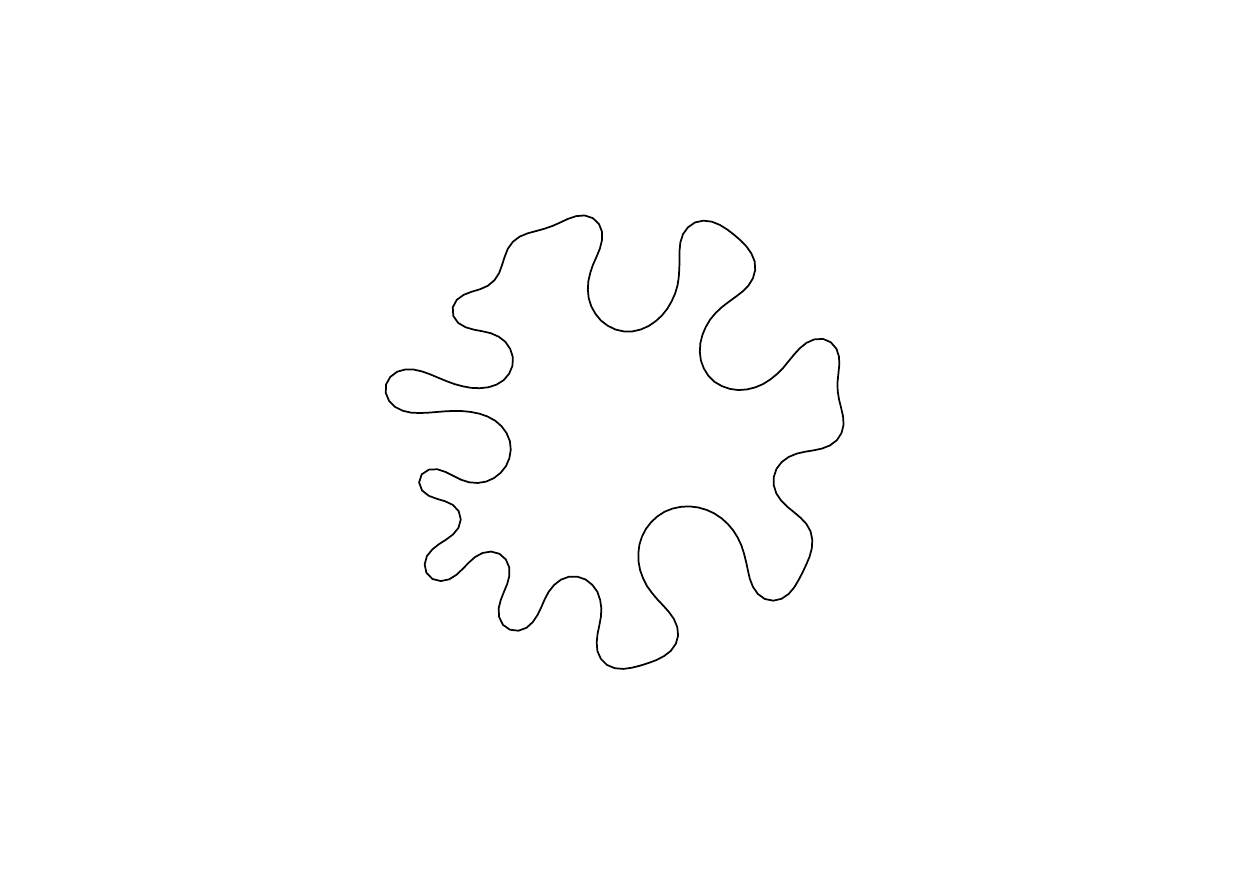}
                \subcaption{$t=0.84$}
            \end{minipage}
            \begin{minipage}{0.25\linewidth}
                \centering
                \includegraphics[scale = 0.32]{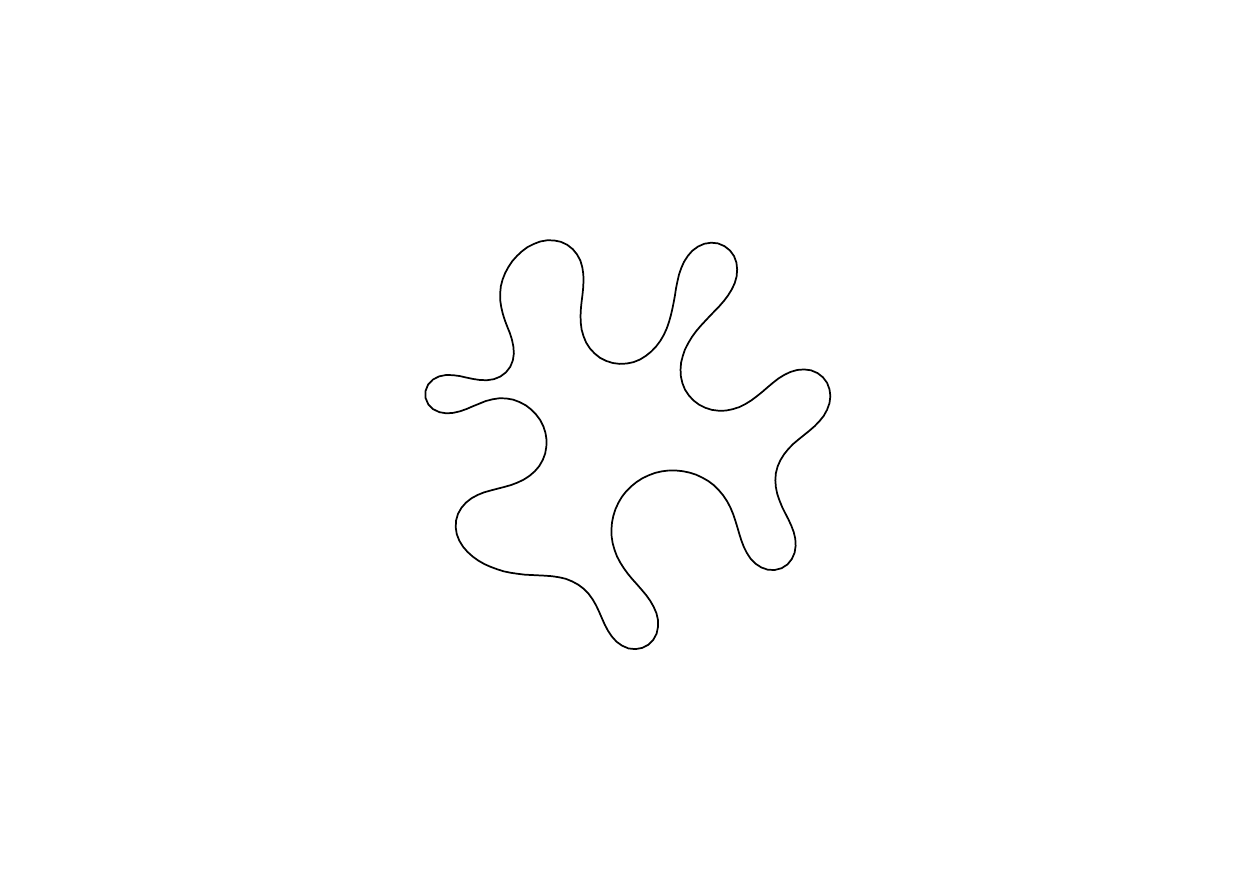}
                \subcaption{$t=1.26$}
            \end{minipage}
            \\
            \begin{minipage}{0.25\linewidth}
                \centering
                \includegraphics[scale = 0.32]{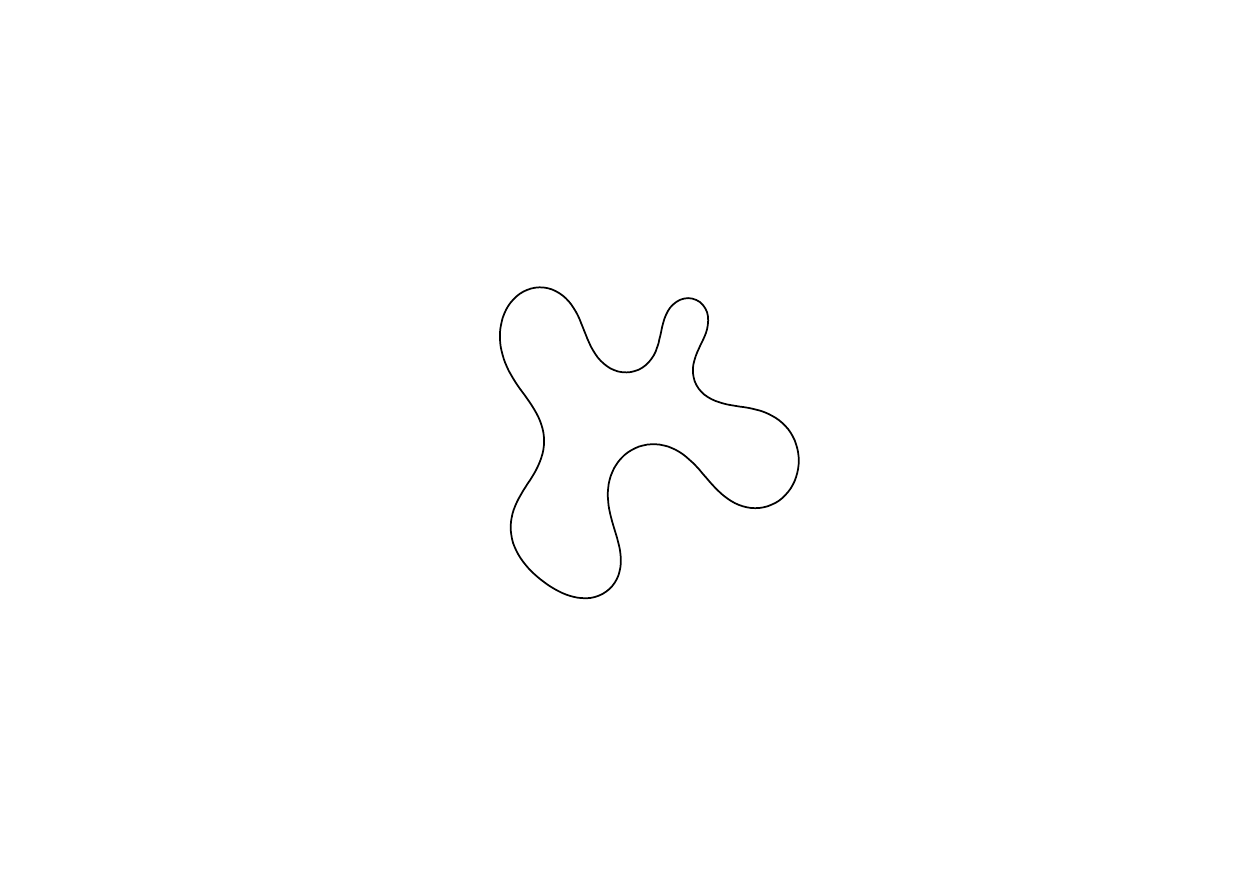}
                \subcaption{$t=1.68$}
            \end{minipage}
            \begin{minipage}{0.25\linewidth}
                \centering
                \includegraphics[scale = 0.32]{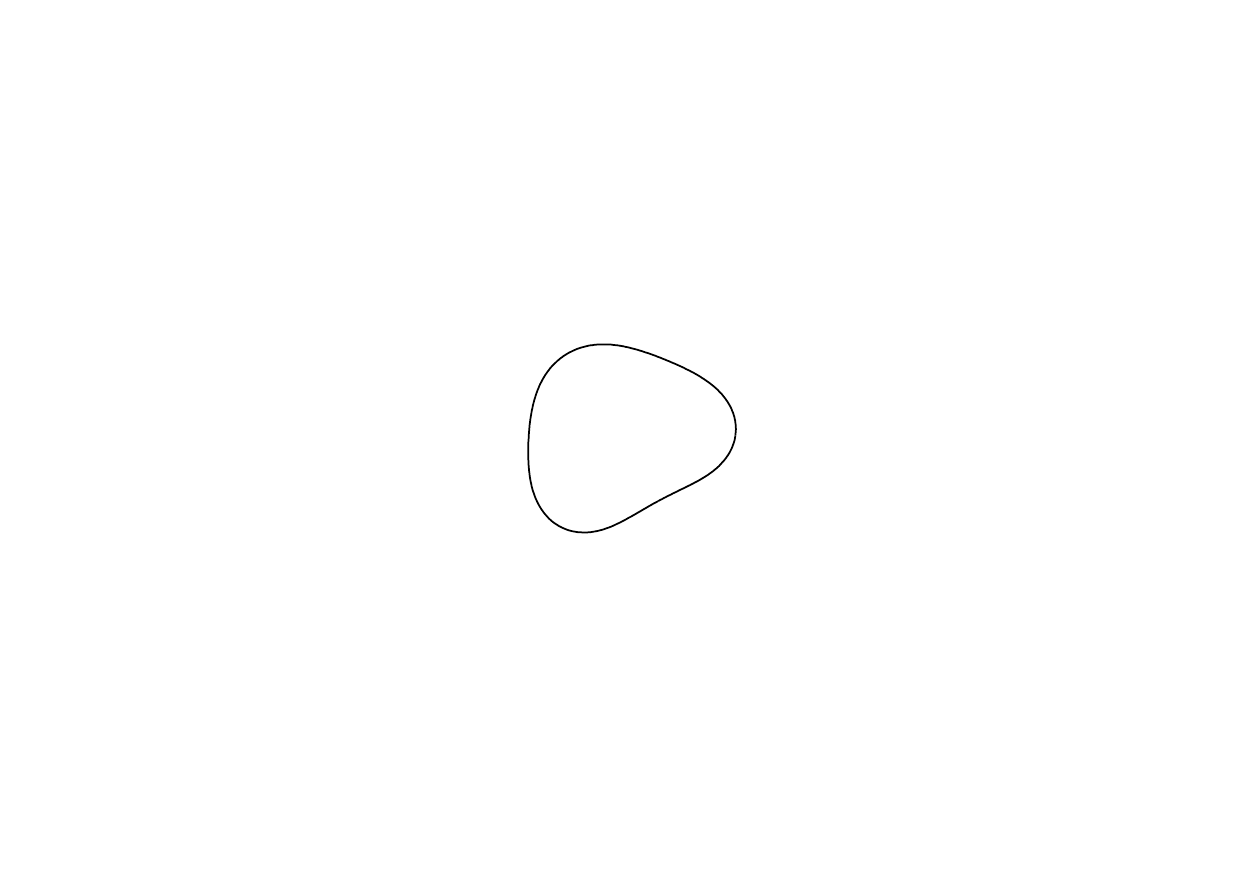}
                \subcaption{$t=2.10$}
            \end{minipage}
            \begin{minipage}{0.25\linewidth}
                \centering
                \includegraphics[scale = 0.32]{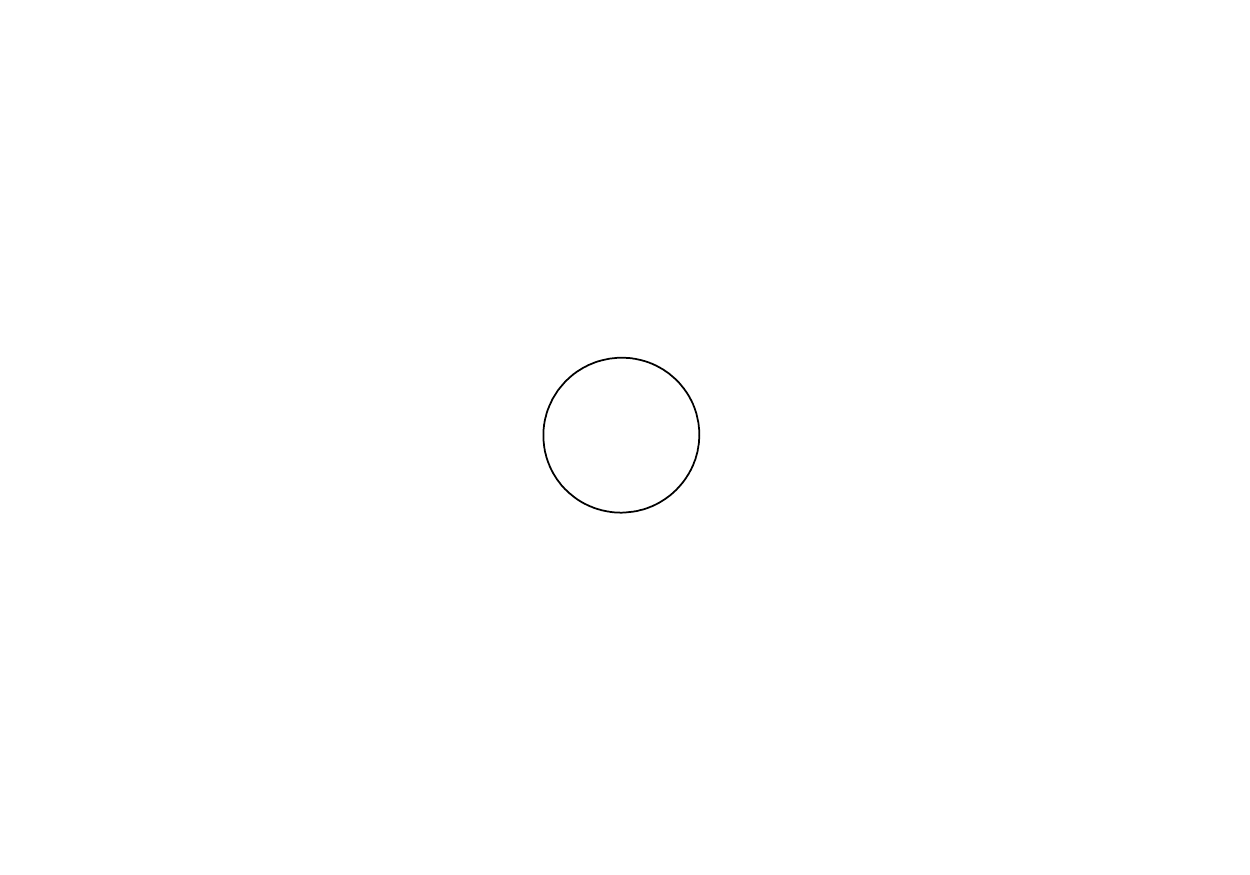}
                \subcaption{$t=2.52$}
            \end{minipage}
            \begin{minipage}{0.25\linewidth}
                \centering
                \includegraphics[scale = 0.32]{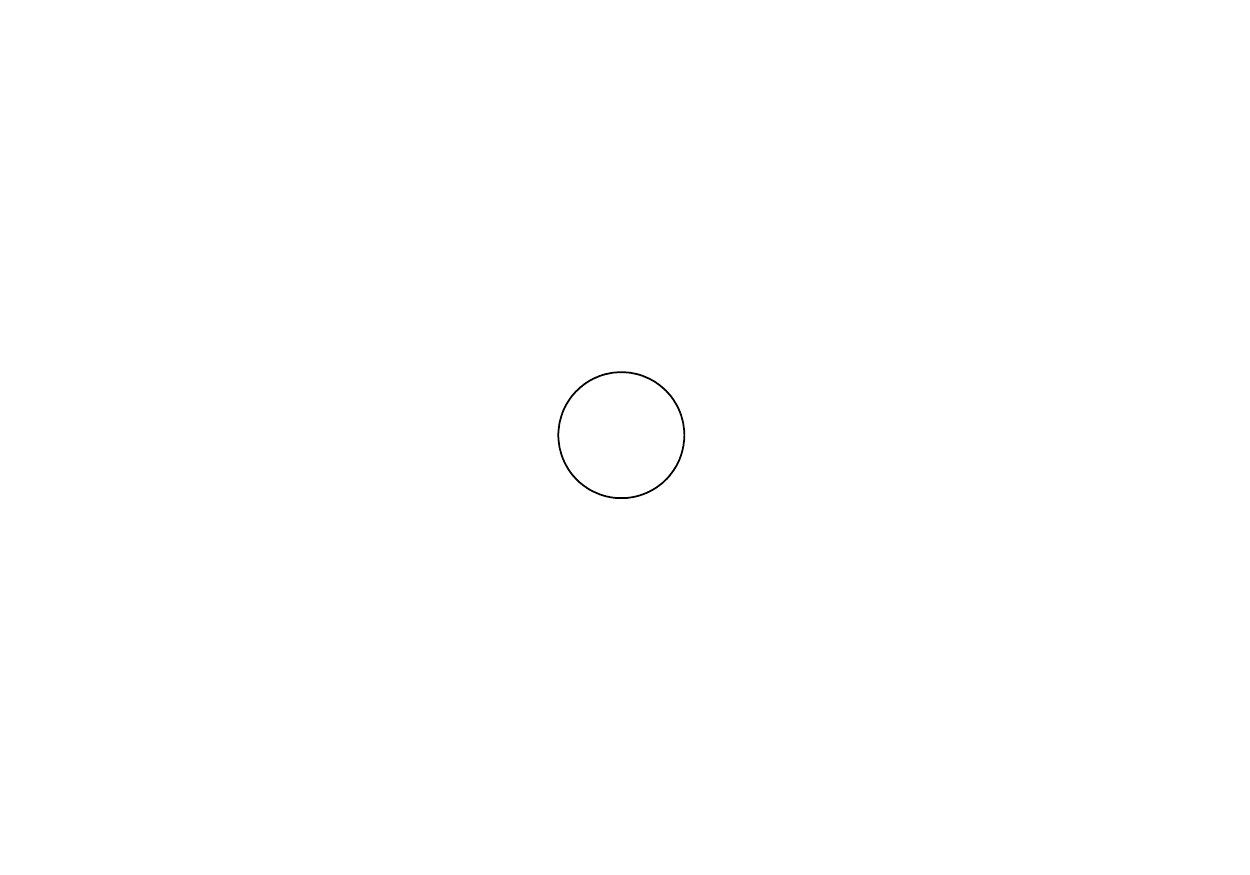}
                \subcaption{$t=2.94$}
            \end{minipage}
        \end{tabular}
            \caption{\label{fig:Bmv0Ca50}$\mathrm{Bmv} = 0, \mathrm{Ca} = 50$}
    \end{figure}

\begin{figure}[tb]
        \begin{tabular}{cccc}
            \begin{minipage}{0.25\linewidth}
                \centering
                \includegraphics[scale = 0.32]{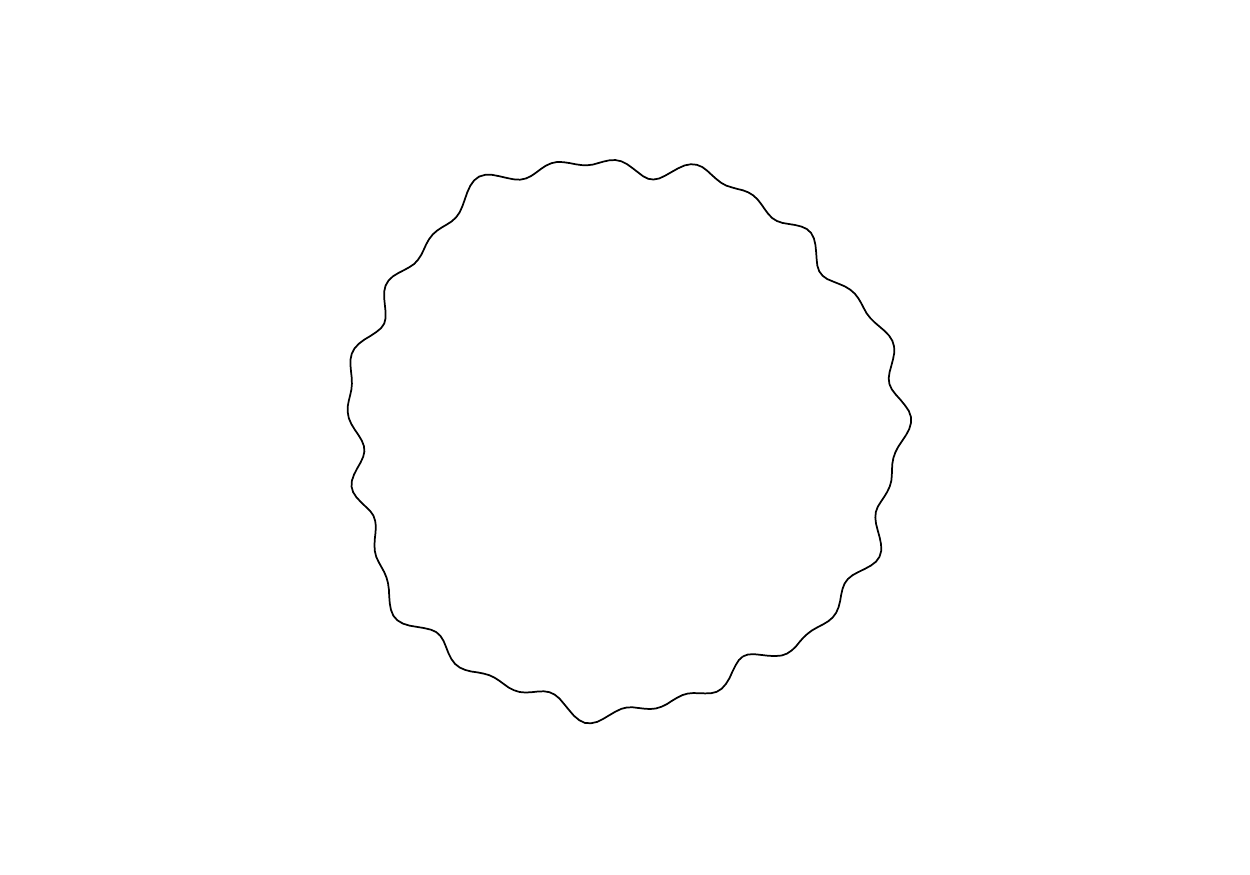}
                \subcaption{$t=0.0$}
            \end{minipage}
            \begin{minipage}{0.25\linewidth}
                \centering
                \includegraphics[scale = 0.32]{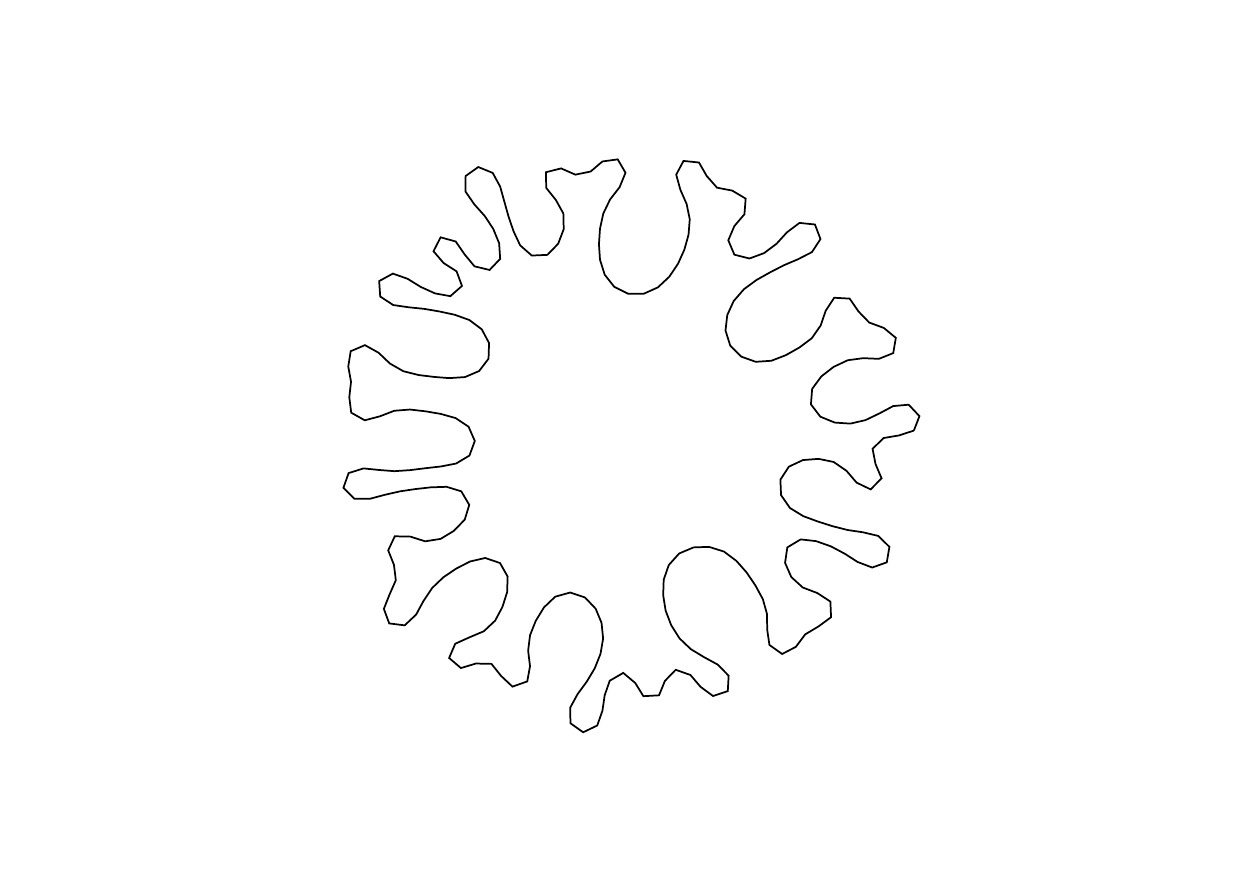}
                \subcaption{$t=0.42$}
            \end{minipage}
            \begin{minipage}{0.25\linewidth}
                \centering
                \includegraphics[scale = 0.32]{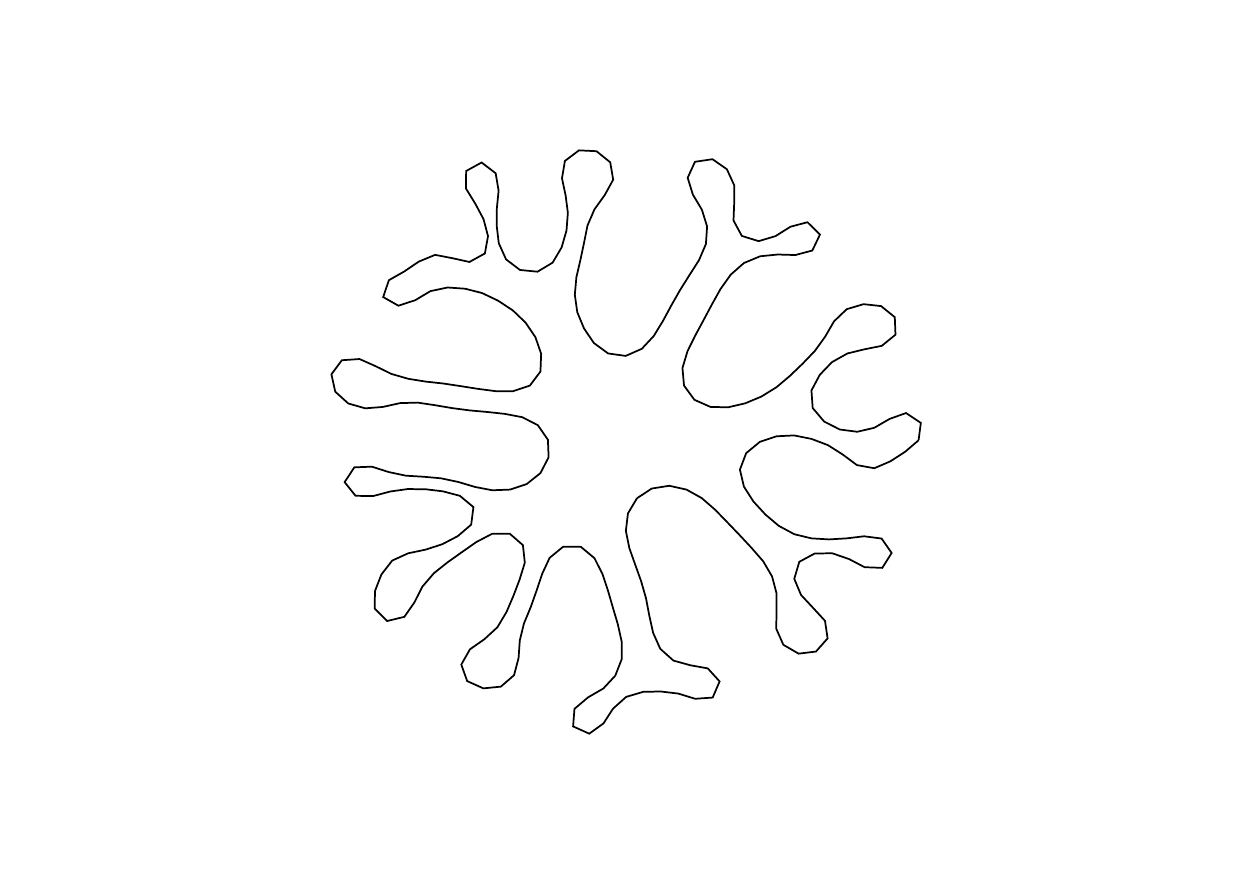}
                \subcaption{$t=0.84$}
            \end{minipage}
            \begin{minipage}{0.25\linewidth}
                \centering
                \includegraphics[scale = 0.32]{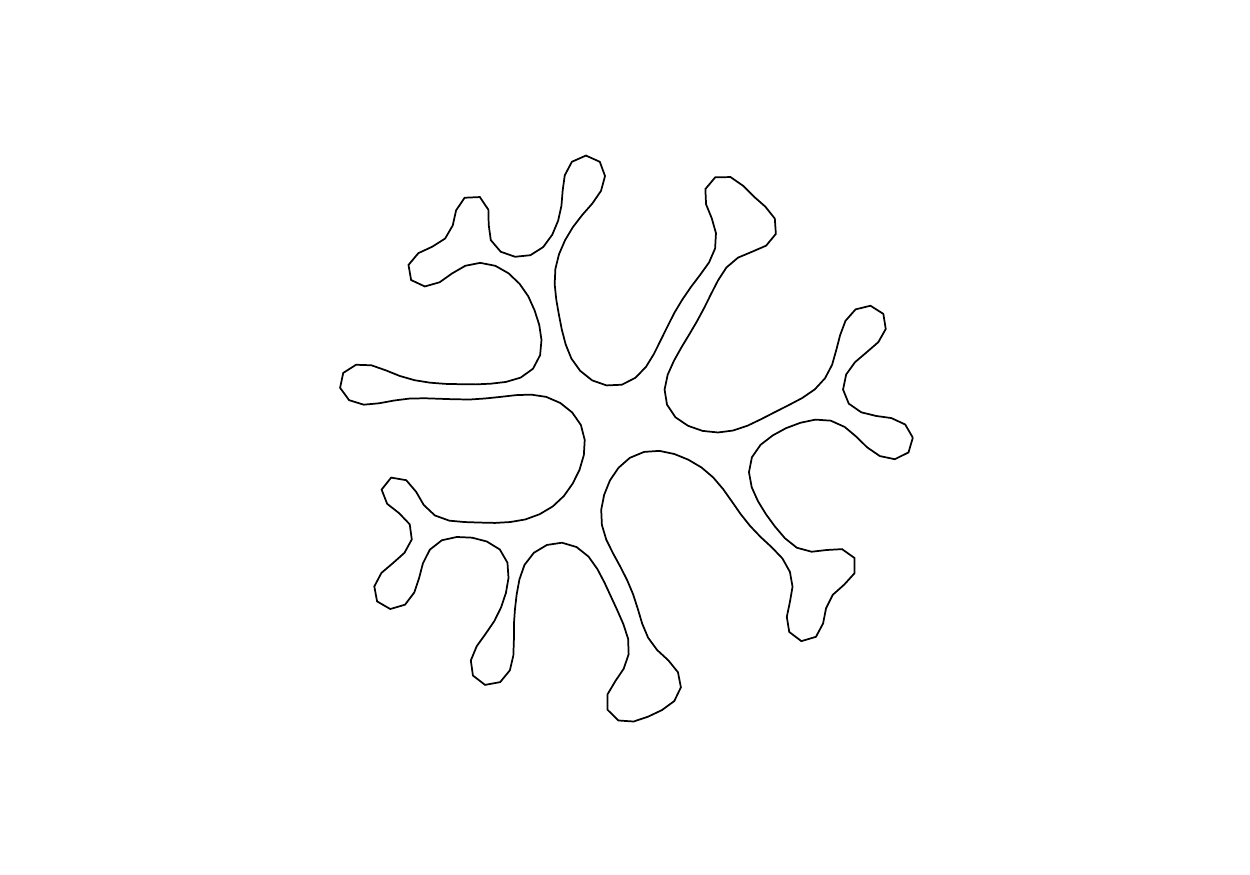}
                \subcaption{$t=1.26$}
            \end{minipage}
            \\
            \begin{minipage}{0.25\linewidth}
                \centering
                \includegraphics[scale = 0.32]{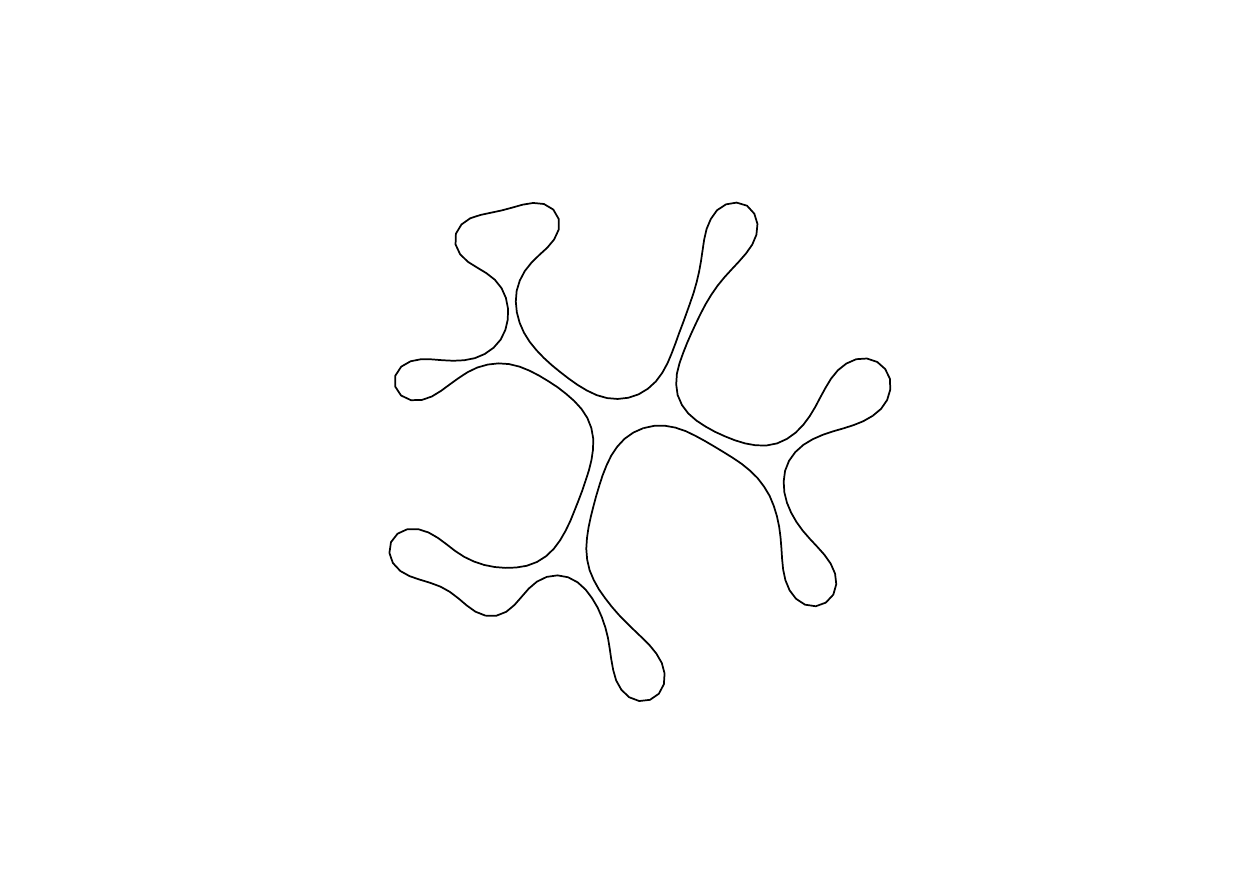}
                \subcaption{$t=1.68$}
            \end{minipage}
            \begin{minipage}{0.25\linewidth}
                \centering
                \includegraphics[scale = 0.32]{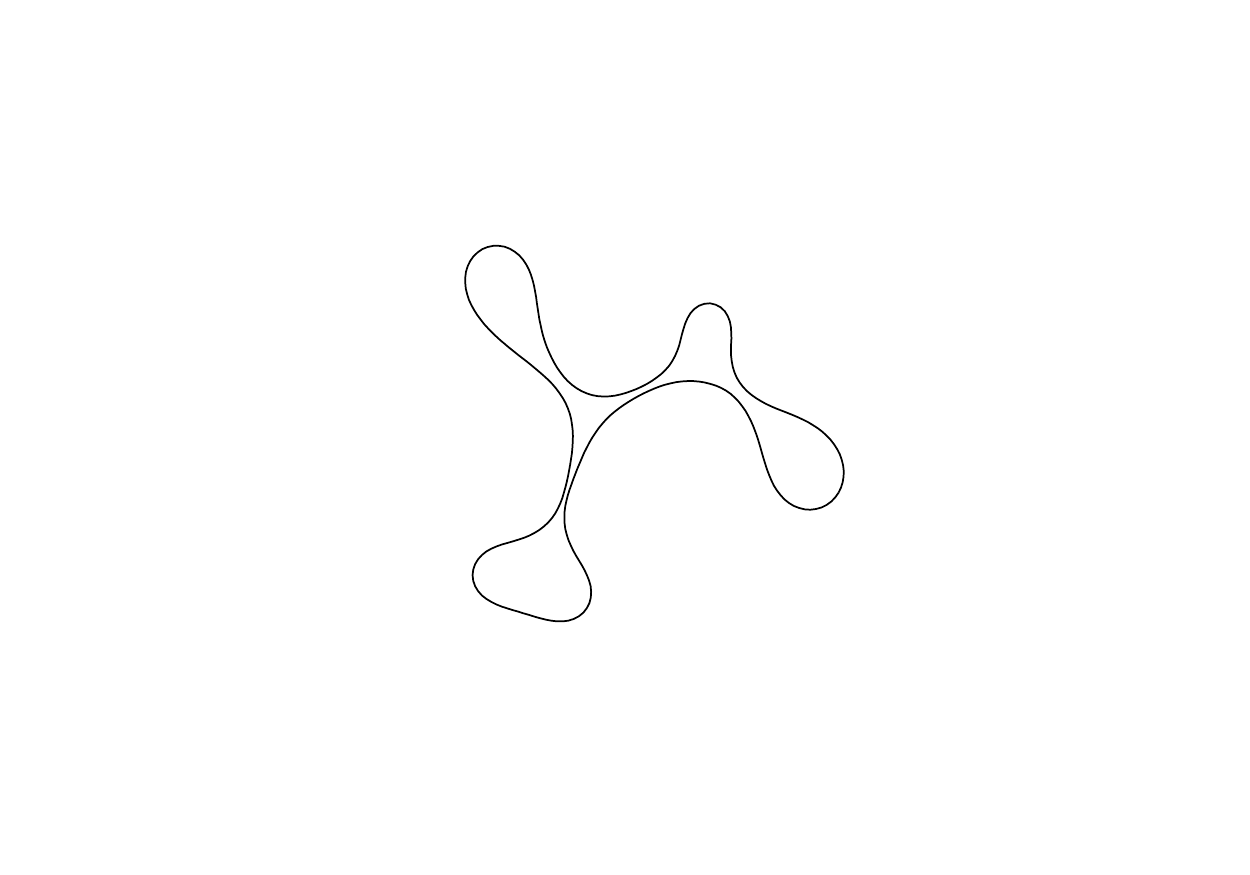}
                \subcaption{$t=2.10$}
            \end{minipage}
            \begin{minipage}{0.25\linewidth}
                \centering
                \includegraphics[scale = 0.32]{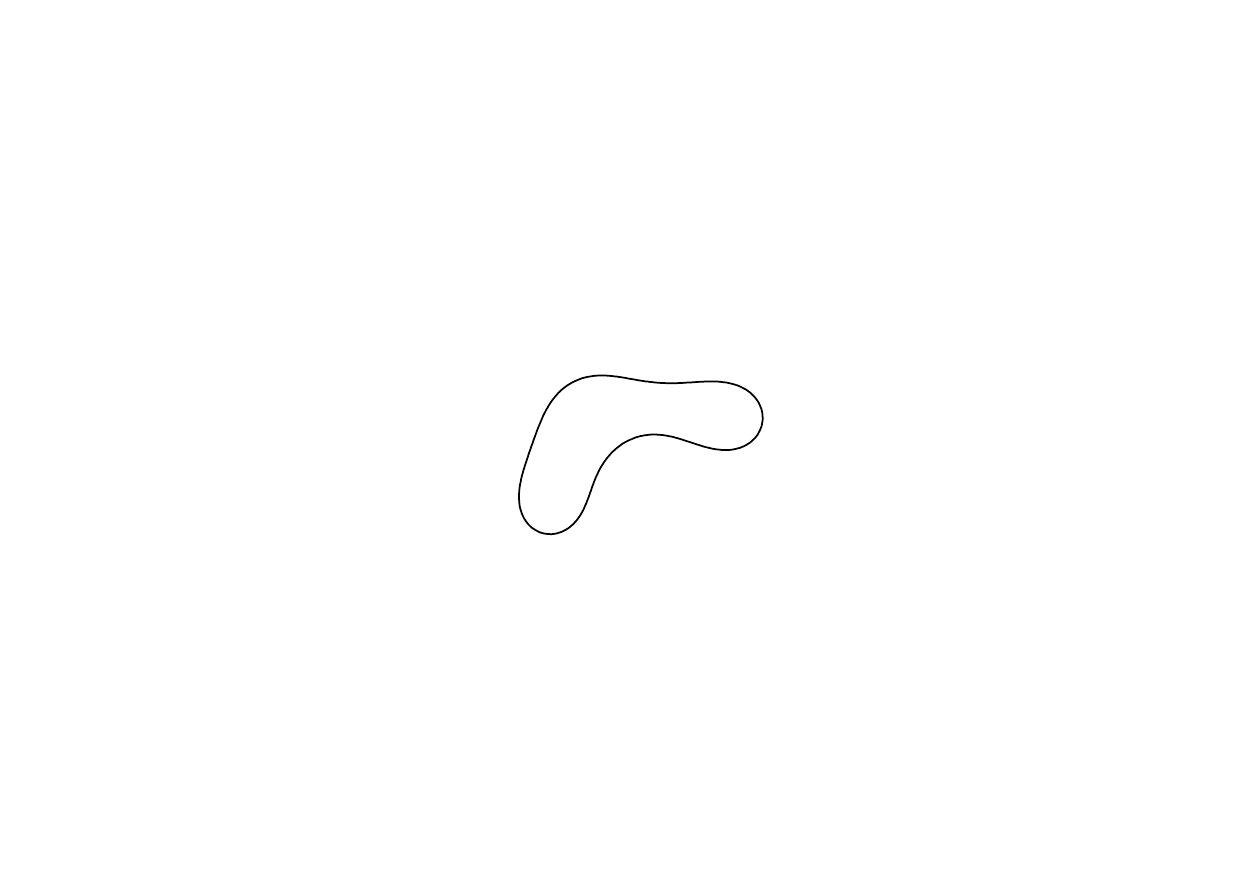}
                \subcaption{$t=2.52$}
            \end{minipage}
            \begin{minipage}{0.25\linewidth}
                \centering
                \includegraphics[scale = 0.32]{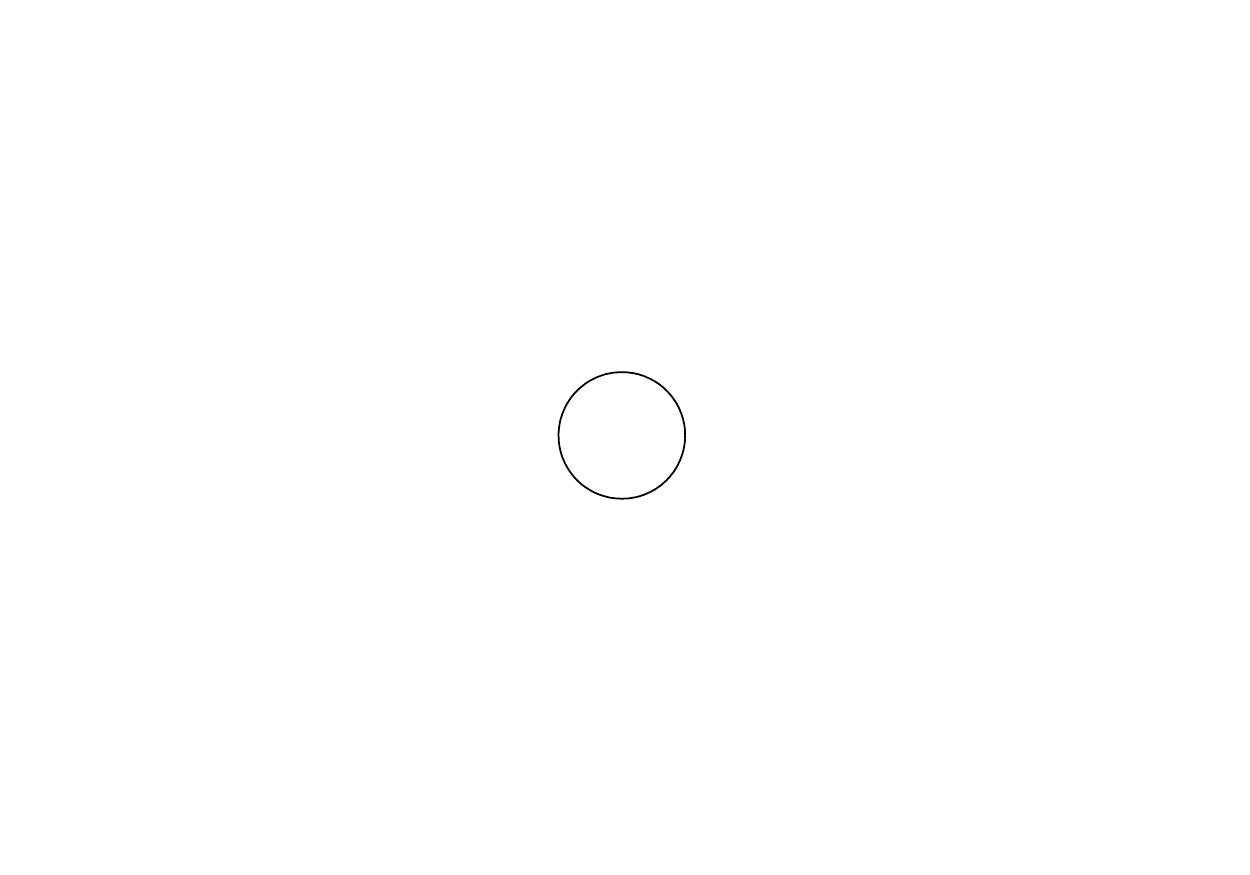}
                \subcaption{$t=2.94$}
            \end{minipage}
        \end{tabular}
        \caption{\label{fig:Bmv35Ca50}$\mathrm{Bmv} = 35, \mathrm{Ca} = 50$}
    \end{figure}

    We show the results of numerical computations for \eqref{eq:HS_magnetic}.
    The initial curve $\Gamma(0) : [0, 1] \ni u \mapsto (x_1(u), x_2(u))^\top \in \mathbb{R}^2$ is given by
    \begin{equation*}
        x_1(u) = r\cos(2\pi u), \quad x_2 = r\sin(2\pi u).
    \end{equation*}
    In Figs. \ref{fig:Bmv0Ca100} and \ref{fig:Bmv25Ca100}, $r$ is given by
    \begin{equation*}
        r = R_0 + 0.02(\cos(6\pi u) + \sin(14\pi u) + \cos(30\pi u) + \sin(50\pi u)),
    \end{equation*}
    and in Figs. \ref{fig:Bmv0Ca50} and \ref{fig:Bmv35Ca50}, $r$ is given 
    \begin{equation*}
        r = R_0 + 0.05(\cos(4\pi u) - \cos(10\pi u) + \cos(22\pi u) - \sin(6\pi u) + \sin(10\pi u)).
    \end{equation*}
    We set the parameters as in Table~\ref{tab_MF}.
    \begin{table}[tb]
        % \centering
        \caption{\label{tab_MF}parameters for the Hele-Shaw flow for magnetic fluids \eqref{eq:HS_magnetic}}
        \begin{tabular}{l|l}
        \hline \hline
        $N = 300$ & the number of vertices\\
        $\Delta t = 1 / (10N^2)$ & the time-mesh size\\
        $r_a = 1$ & the parameter in the Amano's method \eqref{eq:singular_Amano}\\
        $R_0 = 1$ & the initial reference radius\\
        $h_r = h_0 / R_0 = 0.25$ & the parameter introduced for simplifying the model\\
        $M = 1000$ & the number of sample points in the Monte Carlo method\\
        \hline
        \end{tabular}
    \end{table}
    
    In each Fig. \ref{fig:Bmv0Ca100}, \ref{fig:Bmv25Ca100}, \ref{fig:Bmv0Ca50}, and \ref{fig:Bmv35Ca50}, the time evolution of the boundary $\Gamma(t)$ is shown at $t = 0$, $0.42$, $0.84$, $1.26$, $1.68$, $2.10$, $2.52$, and $2.94$.
    Comparing cases $\mathrm{Bmv}=0$, $25$, and $35$, where $\mathrm{Bmv}=(R_0/h_0)^{1/3} \Bm$, it can be seen that the complexity of the pattern formed by the magnetic fluid changes depending on the strength of the magnetic fluid.
    Moreover, the shapes at each time almost coincide with earlier studies by Tatulchenkov and Cebers \cite{tatulchenkov}.
    Note that these are not comparisons at the same discrete time because we use a different time discretization than in previous studies.
    
\begin{figure}[tb]
        \centering
        \includegraphics[scale = 0.5]{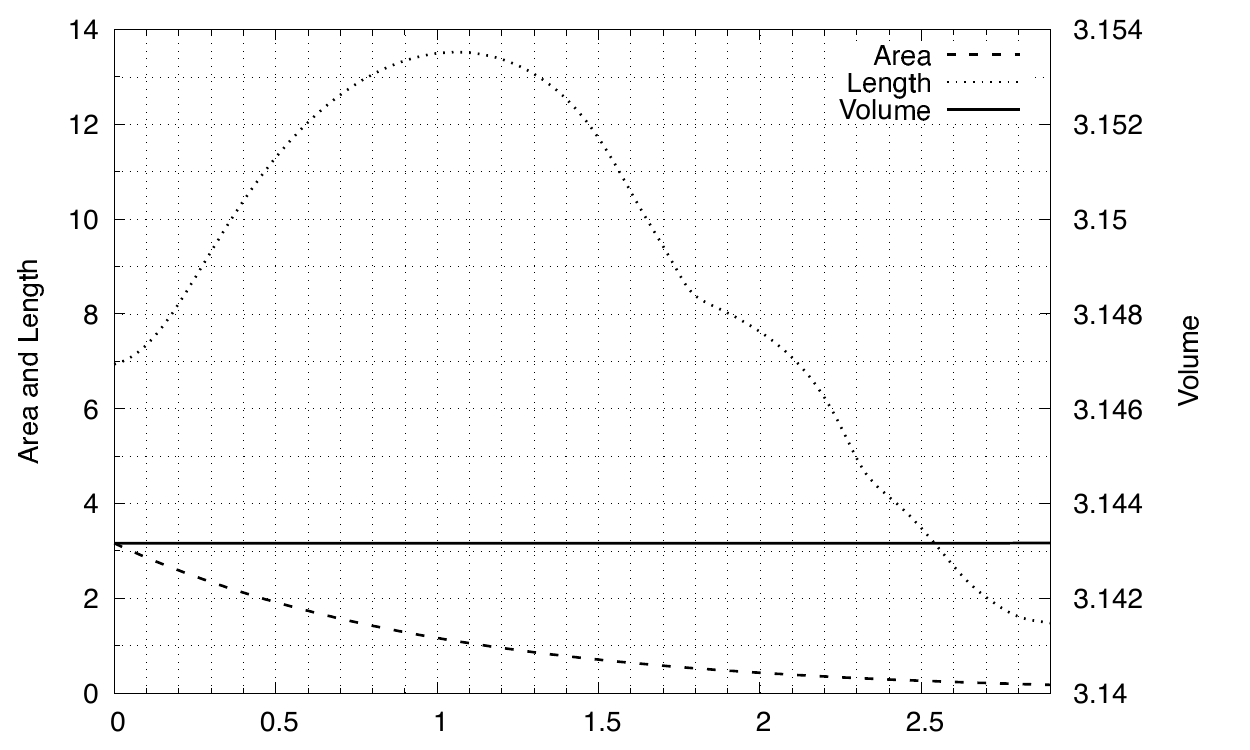}
        \caption{\label{fig:Bmv25Ca100_vol_preserving}time evolution of the Area, the total Length and the Volume ($\mathrm{Bmv} = 25, \mathrm{Ca} = 100$)}    
    \end{figure}

\begin{figure}[tb]
        \centering
        \includegraphics[scale = 0.5]{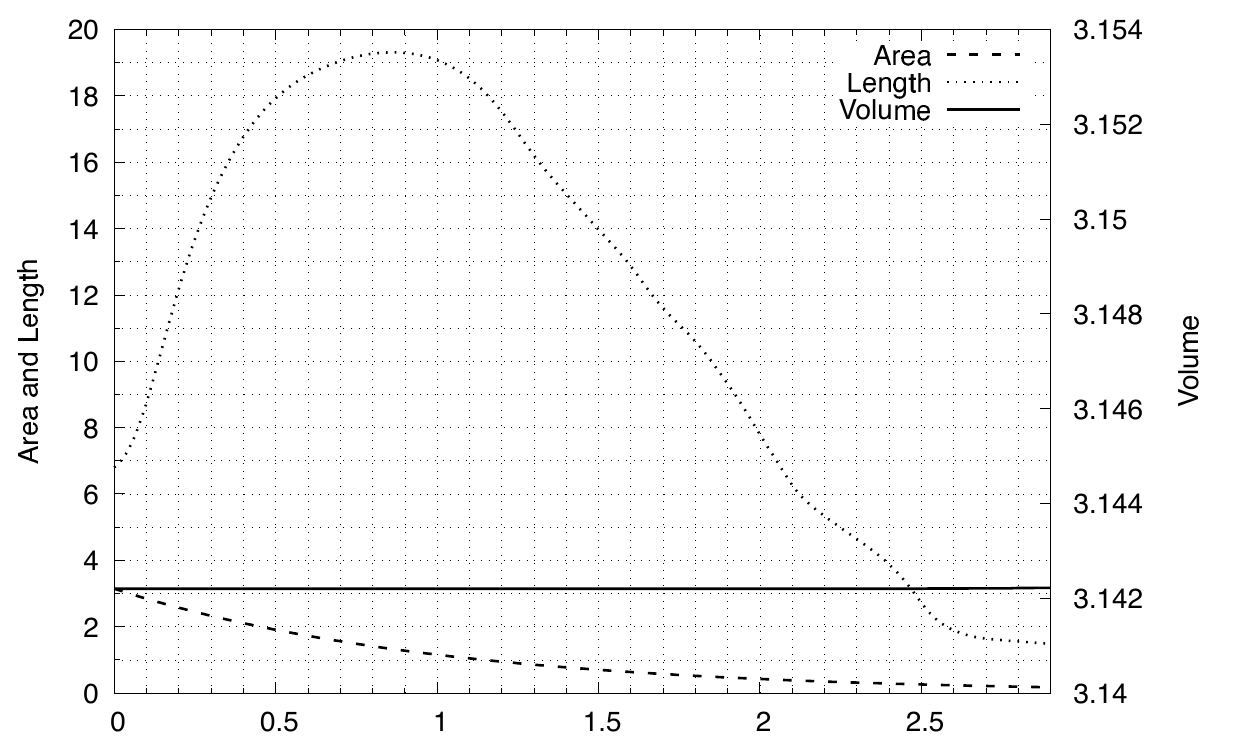}
        \caption{\label{fig:Bmv35Ca50_vol_preserving}time evolution of  the Area, the total Length and the Volume ($\mathrm{Bmv} = 35, \mathrm{Ca} = 50$)}    
    \end{figure}

    Figs. \ref{fig:Bmv25Ca100_vol_preserving} and \ref{fig:Bmv35Ca50_vol_preserving} show the time evolution of the total length, the area, and the volume, where the horizontal axis represents the time, the left vertical axis represents the total length and the area, and the right vertical axis represents the volume.
    The graph shows that the volume-preserving property holds.
    
    
    
    \section{\label{conclusion}Conclusion}
    In this paper, several Hele-Shaw problems were solved using the MFS combined with Amano's method and the asymptotic uniform distribution method.
    The shown results coincided with preceding studies; moreover, some of the properties of Hele-Shaw flow, for example, volume-preserving property, are precisely satisfied.
    Therefore, the results ensured that our numerical method can be adequately effective for Hele-Shaw problems.
    Application of the MFS in moving boundary problems has been reported in quite a few papers; however, our present attempt will widen the possibility of the MFS.
    It is expected that the method will be applied to Hele-Shaw type problems and any other types of potential problems. 
    
    Besides, magnetic fluid instabilities in the Hele-Shaw cell were simulated, then the magnetic effects were indeed observed by comparing numerical results.
    For large $\Bm$, the number of fingers of fingering patterns increased. It is similar to the phenomena in three-dimensional space; that is,  the number of spikes increases with the high intensity of the magnetic field.
    On the other hand, complex patterns such as our results do not appear in three-dimensional space.
    To investigate the difference in the behavior of magnetic fluids between two-dimensional and three-dimensional spaces mathematically is interesting.
    We hope that new aspects for investigations into magnetic fluid are found from our results.

    \bibliographystyle{plain}
    \bibliography{hele-shaw.bib}
\end{document}